\documentclass[11pt]{article}
 
\usepackage{amssymb}
\usepackage{amsmath}
\usepackage{epsfig}
\usepackage{subfigure}
\usepackage{cite}

\newcommand{\cP}{\mathcal P}
\newcommand{\cM}{\mathcal M}
\newcommand{\cA}{\mathcal A}
\newcommand{\cF}{\mathcal F}
\newcommand{\cH}{\mathcal H}
\newcommand{\cO}{\mathcal O}
\newcommand{\cQ}{\mathcal Q}
\newcommand{\cB}{\mathcal B}
\newcommand{\cD}{\mathcal D}
\newcommand{\cR}{\mathcal R}
\newcommand{\cE}{\mathcal E}
\newcommand{\R}{\mathbb{R}}
\newcommand{\Z}{\mathbb{Z}}
\newcommand{\C}{\mathbb{C}} 
\newcommand{\N}{\mathbb{N}}
\newcommand{\de}{\partial}
\newcommand{\itg}{\int \limits}
\newcommand{\e}{\varepsilon}
\newcommand{\ee}{{\rm e}}
\newcommand{\dd}{\sqrt{\cD}}
\newcommand{\ds}{\displaystyle}

\DeclareMathOperator{\dist}{dist}
       
\newtheorem{theorem}{Theorem}[section]
\newtheorem{example}{Example}[section]
\newtheorem{remark}{Remark}[section]
\newtheorem{corollary}{Corollary}[section]

\begin{document}

\markboth{F. Lanzara et al.}{Approximate Hermite quasi-interpolation}

\title{Approximate Hermite quasi-interpolation}

\author{{  F. Lanzara$^{\mbox{\tiny 1}}$ , V. Maz'ya$^{\mbox{\tiny 2}}$ ,
G. Schmidt$^{\mbox{\tiny 3}}$}
}
\date{}
\maketitle
\baselineskip=0.9\normalbaselineskip
\vspace{-3pt}

\hspace*{1cm}
\parbox{10cm}{\begin{flushleft}
{\footnotesize\em
\begin{itemize}
\item[$^{\mbox{\tiny\rm 1}}$]Mathematics Department, Sapienza
University of Rome,\\
Piazzale Aldo Moro 2, 00185 Rome, Italy\\
\texttt{\rm lanzara\symbol{'100}mat.uniroma1.it}
\item[$^{\mbox{\tiny\rm 2}}$]Department of Mathematics, University of
Link\"oping, \\ 581 83 Link\"oping, Sweden;\\
Department of Mathematics, Ohio State University,\\
231 W 18th Avenue, Columbus, OH 43210, USA;\\
Department of Mathematical Sciences, M\&O Building, University of
Liverpool, Liverpool L69 3BX, UK;\\
\texttt{\rm vlmaz\symbol{'100}mai.liu.se }
\item[$^{\mbox{\tiny\rm 3}}$]Weierstrass Institute for Applied Analysis and
Stochastics, \\  Mohrenstr. 39,
10117 Berlin, Germany \\
\texttt{\rm schmidt\symbol{'100}wias-berlin.de}
\end{itemize}
}
\end{flushleft}}

{\bf Abstract.}
In this paper we derive approximate quasi-interpolants when
the values of a function $u$ and of some of its derivatives are prescribed
at the points of a uniform grid. As a byproduct of these formulas we
obtain very simple approximants  which provide
high order approximations for solutions to elliptic differential 
equations with constant coefficients.

\bigskip


\bigskip


\section{Introduction}\label{S:intro}
 
The method of approximate approximations 
is mainly directed to the numerical solution of partial
integro-differential equations. 
The method provides simple formulas for
quasi-interpolants, which approximate functions up to a prescribed precision very accurately,
but in general the approximants do not converge. The lack of
convergence, which is not 
perceptible in numerical computations, is offset by a greater flexibility in the choice of approximating
functions. So it is possible to construct multivariate approximation formulas, which are easy
to implement and have additionally the property that
pseudodifferential operations can be 
effectively performed. This allows to create effective numerical algorithms for solving boundary
value problems for differential and integral equations.

Approximate quasi-interpolants on a uniform grid are of the form
\begin{equation}\label{general0}
 Mu(x)=\cD^{-n/2}\sum_{j\in\Z^{n}}\cH\left(\frac{x-h j}{h\sqrt{\cD}}
\right)u(hj)
\end{equation}
with positive parameters, ``small'' $h$ and ``large'' $\cD$,  
and the generating
function $\cH$
is sufficiently smooth and of rapid decay.
Their properties have been 
studied in a  
series of papers
\cite{Ma2,MS1,MS2,Ma1} (cf. also the monograph \cite{MSbook}).
It has been shown
that the quasi-interpolant approximates smooth functions with
\begin{equation}\label{generalest}
|Mu(x)-u(x)|\leq
\e \sum_{|\beta|=0}^{N-1} (h
         \sqrt{\cD})^{|\beta|}|\de^{\beta} u(x)|+ c (h \sqrt{\cD})^{N}
         \Vert \nabla_{N} u \Vert_{L_{\infty}} \, ,
\end{equation}
as long as $\cH$ is subject to the moment condition
\begin{equation}\label{momcon}
\itg_{R^n} x^\alpha \cH(x) \, dx = \delta_{|\alpha|0}\, , \quad 0 \le |\alpha| < N \, .
 \end{equation}
In \eqref {generalest} the constant $c$ depends only on
$\cH$ and $\e$ can be made arbitrarily small if the
 parameter $\cD$ is
sufficiently large, so that one can fix  $\cD$ such that in numerical computations
$Mu$ approximates with the order $O(h^N)$.
The construction of simple generating functions satisfying \eqref{momcon}
for arbitrary $N$
has been addressed in \cite{MS3}, whereas  \cite{LMS,MS4}
extend the results to the case of nonuniform grids.

In this paper 
we study more general quasi-interpolation formulas
with values of a function $u$ and of some of its derivatives prescribed
at the points of a uniform grid.
More precisely, we consider approximants of the form 
\begin{equation}\label{generalqi}
Mu(x)=\cD^{-n/2}\sum_{j\in\Z^{n}}\cH\left(\frac{x-h j}{h\sqrt{\cD}}
\right)\cQ((-h\sqrt{\cD})\partial)u(hj)\, ,
\end{equation}
where $\cQ(t),\, t\in \R^{n}$, is a polynomial with $\deg \cQ<N$.

We establish estimates of the type \eqref{generalest}
if the generating function $\cH$ and the coefficients of $\cQ$ are
connected by suitable
conditions (see \eqref{3.10}). 
It is shown in particular, that
for an arbitrary polynomial $\cQ$ there exists $\cH$ such that
the approximate Hermite quasi-interpolant \eqref{generalqi}
satisfies the estimate \eqref{generalest}.
On the other hand, 
the same is true
if $\cH$ is  the Gaussian or a related function
and $\cQ$ is chosen suitably.
As a byproduct
we obtain very simple
approximants  which provide
high order approximations to solutions of the Laplace equation or
other elliptic differential equations.

Estimate \eqref{generalest} is proved in Section \ref{S:DUE} 
for  \eqref{generalqi} with
$\cH$ and  $\cQ$
connected by 
conditions \eqref{3.10}.
In Section \ref{S:QUA} we describe a simple method for the
construction of $\cH$ satisfying the conditions \eqref{3.10} for a
given polynomial $\cQ$. This result is applied in
  Section \ref{example} when we consider some examples of Hermite 
quasi-interpolants.
In the particular case $\cH(x)=\pi^{-n/2}(\det B)^{-1/2}{\rm
e}^{-<B^{-1}x,x>}$, with $B=\{b_{ij}\}$ symmetric and positive definite
real matrix, we obtain
the following quasi-interpolant of order
$\cO((h\dd)^{2M})$
\[M u(x)=\frac{(\det B)^{-1/2}}{(\pi \cD)^{n/2}}
\sum_{m\in\Z^{n}}\sum_{s=0}^{M-1}
\frac{(-h^2 \cD )^{s}}{s!4^{s}} \, \cB ^{s}u(hm)\,\ee^{\, -\langle
   B^{-1}(x-hm),\,x-hm\rangle/(h^2\cD)}
\]
where $\cB$ is the second order partial differential operator
$\cB u=\sum_{i,j=1}^{n}b_{ij}\de_{i}\de_{j}u$.
In Section \ref{appl} we use $M u$
for the approximation of solutions of the equation $\cB u=0$.
If $u$ satisfies the equation in $\R^{n}$
then  $M u$  has the simple form of the quasi-interpolant
of order $\cO((h\dd)^{2})$
but  gives an approximation of exponential order plus 
a small saturation error.
The same approximation property holds  for functions which satisfy 
the equation $\cB u=0$ in a domain $\Omega\subset \R^{n}$. If the 
function $u$ is
extended by zero outside $\Omega$ then $M u$ simplifies to
\[
M u(x)=\frac{(\det B)^{-1/2}}{(\pi \cD)^{n/2}}
\sum_{hm\in\Omega}\, u(hm)\,\ee^{\, -\langle
   B^{-1}(x-hm),\,x-hm\rangle/(h^2\cD)} .
\]
We obtain that, for any $\e>0$, $M\,u$ approximates pointwise $u$ in 
a subdomain
$\Omega'\subsetneq \Omega$
with
$$|M u(x)-u(x)|\leq
\e \sum_{|\beta|=0}^{2M-1} (h
         \sqrt{\cD})^{|\beta|}|\de^{\beta} u(x)|+ C_B (h \sqrt{\cD})^{2M}
        c_{u},$$
if we choose $\cD$ and $h$ appropriately; the constant $C_B$ is independent of
$u,h,\cD, M$ and $c_{u}$ depends on $u$.

In Section \ref{deriv} we show that the Hermite quasi-interpolant
\eqref{generalqi} gives the simultaneous
approximation of the derivatives of $u$. If $\de^{\beta}\cH$ exists and 
is of rapid decay,
for any function $u\in W_{\infty}^{L}(\R^{n})$ with $L\geq N+|\beta|$,
the difference $\de^{\beta}Mu(x)-\de^{\beta}u(x)$
can be estimated by
\[
|\de^{\beta}Mu(x)-\de^{\beta}u(x)|\leq
  \e\,\sum_{|\gamma|=0}^{L-1}(h\dd)^{|\gamma|-|\beta|}|\de^{\gamma}u(x)|+
  c_{1}\,(h\dd)^{N}\Vert \nabla_{N+|\beta|} u \Vert_{L_{\infty}},
\]
with $c_{1}$
independent of $u$, $h$ and $\cD$.

\section{Quasi-interpolants with derivatives}\label{S:DUE}

In this section we study the approximation of a smooth function
$u(x)$, $x\in \R^{n}$,
by the  Hermite quasi-interpolation operator
\begin{equation}\label{3.0}
       M u(x)=\cD^{-n/2} \sum_{m\in \Z^{n}}\cH\left( \frac{x-h m}{h
\sqrt{\cD}}\right)
           \cQ\left(-h
           \sqrt{\cD}\, {\partial }\right)  u(hm)
\end{equation}
where $\cQ(t)$, $t\in \R^{n}$, is a polynomial of degree at most $N-1$
\begin{equation}\label{polinomio}
\cQ(t)=\sum_{|\gamma|=0}^{N-1}a_{\gamma} t^{\gamma},\> t \in \R^{n},
\quad \mbox{with }a_{\gamma}\in \R \; \mbox{ and }
a_{0}=1 \, ,
\end{equation}
$\partial=\partial_1 \ldots\partial_n$,
and $\cH$ is a sufficiently smooth, rapidly decaying function.
Then
\begin{equation}
       M u(x)=\cD^{-n/2} \sum_{m\in \Z^{n}}\left(\sum_{|\gamma|=0}^{N-1}(-h
       \sqrt{\cD})^{|\gamma|} a_{\gamma} \de^{\gamma} u(hm) \right)
\cH\left( \frac{x-h m}{h \sqrt{\cD}}
           \right).
       \label{3.1}
\end{equation}

Our aim is to give conditions on the generating function
$\cH$ such that (\ref{3.1}) is an
approximation formula of order $\cO((h \sqrt{\cD})^{N})$ plus
terms, which can be made sufficiently small.

Suppose that $u\in C^{N}(\R^{n})\cap L^\infty(\R^n)$.
Taking the Taylor expansion of
$\de^{\gamma} u$ at each node $h m$ leads to
    \begin{equation} \label{3.2}
    \begin{split}
        \de^{\gamma}u(h m)= &\sum_{|\alpha|=0}^{N-1-|\gamma|} \frac{(h m-x)^{\alpha}}{\alpha!} \,
        \de^{\alpha+\gamma} u(x)  \\
    &+\sum_{|\alpha|=N-|\gamma|} \frac{(h
        m-x)^{\alpha}}{\alpha!} \, U_{\alpha+\gamma}(x,h m),\quad 0\leq |\gamma|<N,
    \end{split}
    \end{equation}
with
\begin{equation}\label{ualfa}
       U_{\alpha}(x,y)=N \itg_{0}^{1}s^{N-1}\partial^{\alpha} u(s x+(1-s)y)\,
       ds\,, \;
|\alpha|=N.
\end{equation}
We write $Mu$ in the form
\begin{equation*}
       M u(x)= \sum_{|\beta|=0}^{N-1} (- h \sqrt{\cD})^{|\beta|}\de^{\beta}u(x)
       \sum_{\alpha\leq \beta}
       \frac{ a_{\beta-\alpha}}{\alpha!}
       \sigma_{\alpha}(\frac{x}{h},\cD,\cH)+\cR_{h,N}(x)
\end{equation*}
with the periodic functions
\begin{equation*}
       \sigma_{\alpha}(\xi,\cD,\cH)=\cD^{-n/2} \sum_{m \in \Z^{n}} \left( \frac{\xi-
       m}{\sqrt{\cD}}\right)^{\alpha} \cH\left( \frac{\xi-
       m}{\sqrt{\cD}}\right),\quad 0\leq |\alpha|<N \, ,
\end{equation*}
and the remainder term
\begin{equation*}
\begin{split}
&\cR_{h,N}(x)= \\
&\; (-h \sqrt{\cD})^{N}
\cD^{-n/2}\sum_{|\beta|=N} \sum_{0<\alpha\leq\beta} \frac{a_{\beta-\alpha}}{\alpha!} \!\!
    \sum_{m\in \Z^{n}} U_{\beta}(x,h m)
\left( \frac{x-h m}{h \sqrt{\cD}}\right)^{\alpha}\!\!\!
    \cH\left( \frac{x-h m}{h \sqrt{\cD}}\right).
\end{split}
\end{equation*}
Therefore by using the definition of $\sigma_{\alpha}$,
\begin{equation}\label{3.6}
\begin{split}
       M u(x)-u(x)=&\,u(x) \left( 
\sigma_{0}\big(\frac{x}{h},\cD,\cH\big)-1\right) \\
&+\sum_{|\beta|=1}^{N-1} (-h \sqrt{\cD})^{|\beta|} \de^{\beta}u(x)
       \sum_{\alpha\leq \beta} \frac{a_{\beta-\alpha}}{\alpha!}
       \sigma_{\alpha}\big(\frac{x}{h},\cD,\cH\big)+\cR_{h,N}(x).
\end{split}
\end{equation}
Let us introduce the functions
\[
       \cE_{\alpha}(\xi,\cD,\cH)=\sigma_{\alpha}(\xi,\cD,\cH)-
    \itg_{\R^{n}} x^{\alpha} \cH(x)\, dx \, ,\qquad 0\leq|\alpha|\leq N-1 .
\]
Then we have
\begin{equation}\label{3.8}
       \sigma_{0}(\xi,\cD,\cH)-1=\cE_{0}(\xi,\cD,\cH)+
       \Big(\itg_{\R^{n}} \cH(x)
       dx-1\Big)
\end{equation}
and for $|\beta|=1,\ldots,N-1$
\begin{equation}\label{3.9}
    \sum_{\alpha\leq \beta}
               \frac{a_{\beta-\alpha}}{\alpha!} \sigma_{\alpha}(\xi,\cD,\cH)
=  \sum_{\alpha\leq \beta}
           \frac{a_{\beta-\alpha}}{\alpha!}
           \cE_{\alpha}(\xi,\cD,\cH)
     + \sum_{\alpha \leq \beta}
               \frac{a_{\beta-\alpha}}{\alpha!}
    \itg_{\R^{n}}x^{\alpha} \cH(x) dx.
\end{equation}

\begin{theorem}\label{UNO} Suppose that $\cH$ is differentiable up to the order
of the smallest integer $n_{0}>n/2$, satisfies the
       decay condition: there exist $K>N+n$ and $C_{\beta}>0$ such that
\begin{equation}\label{decaybis}
       |\partial ^{\beta} \cH(x)|\leq C_{\beta}\, (1+|x|)^{-K},\quad
       0\leq|\beta|\leq n_{0}
\end{equation}
      and the conditions
\begin{equation} \label{3.10}
       \begin{cases}
\ds{\qquad \itg_{\R^{n}} \cH(x) dx }& \hskip-5pt = 1 \\[5mm]
       \ds{\; \sum_{\alpha \leq \beta} \frac{a_{\beta-\alpha}}{\alpha!}  \itg_{\R^{n}}x^{\alpha}
           \cH(x) dx} & \hskip-5pt = 0 \,,\quad |\beta|=1,\ldots,N-1 .
\end{cases}
\end{equation}
       Then for any $\e>0$ there exists $\cD>0$ such that, for all $u
       \in W^{N}_{\infty}(\R^{n})\cap C^{N}(\R^{n})$ the approximation
error of the
       quasi-interpolant \eqref{3.1} can be pointwise estimated by
       \begin{equation}\label{3.14}
           |M u(x) -u(x)|\leq \e \sum_{|\beta|=0}^{N-1} (h
           \sqrt{\cD})^{|\beta|}|\de^{\beta} u(x)|+ c (h \sqrt{\cD})^{N}
           \Vert \nabla_{N} u \Vert_{L_{\infty}}
       \end{equation}
           with the constant $c$ not depending on $h$, $u$ and $\cD$.
\end{theorem}
{\bf Proof.} Under conditions (\ref{3.10}) the relations
\eqref{3.8} and \eqref{3.9} simplify to
\begin{equation*}
       \begin{array}{rl}
       \sigma_{0}(\xi,\cD,\cH)-1 &  = \cE_{0}(\xi,\cD,\cH) , \\ \\
\displaystyle{    \sum_{\alpha\leq \beta}
               \frac{a_{\beta-\alpha}}{\alpha!} \sigma_{\alpha}(\xi,\cD,\cH)}  &
 = \displaystyle{\sum_{\alpha\leq \beta}
              \frac{a_{\beta-\alpha}}{\alpha!}
              \cE_{\alpha}(\xi,\cD,\cH)}.
\end{array}
\end{equation*}
Moreover, the assumptions (\ref{decaybis}) on $\cH$ ensure that
(see \cite{MS4})
\[
\{\partial^{\alpha}\cF
\cH(\sqrt{\cD} \nu)\}\in l_{1}(\Z^{n}),\qquad 0\leq |\alpha|<K-n,
\]
\begin{equation*}
       | \cE_{\alpha}(\xi,\cD,\cH)|=\Big|\sigma_{\alpha}(\xi,\cD,\cH)-
    \int_{\R^{n}} x^{\alpha} \cH(x)\, dx\Big|\leq (2\pi)^{-|\alpha|} \hskip-3pt
    \sum_{\nu\in\Z^{n}\setminus 0}\big|\partial^{\alpha}\cF
\cH(\sqrt{\cD} \nu)\big| ,
\end{equation*}
and
\begin{equation*}
\sum_{\nu\in\Z^{n}\setminus 0}|\partial^{\alpha}\cF
\cH(\sqrt{\cD} \nu)| \longrightarrow 0 \;\mbox{ as }\; \cD \to \infty  .
\end{equation*}
Hence
for any $\e>0$ there exists $\cD$ such that the following
estimates are valid
\begin{equation}
       |\sigma_{0}(\xi,\cD,\cH)-1|=|\cE_{0}(\xi,\cD,\cH)|< \e,
       \label{3.11}
\end{equation}
\begin{equation}
      \Big| \sum_{\alpha\leq \beta}
              \frac{a_{\beta-\alpha}}{\alpha!} \sigma_{\alpha}(\xi,\cD,\cH)\Big| \leq
              \sum_{\alpha\leq \beta}
              \frac{|a_{\beta-\alpha}|}{\alpha!}
              |\cE_{\alpha}(\xi,\cD,\cH)| <\e,\quad |\beta|=1,\ldots,N-1.
       \label{3.12}
\end{equation}
The next step is to estimate the remainder term $\cR_{h,N}$.
Since
\begin{equation*}
       |U_{\alpha}(x,y)|=N\left|\itg_{0}^{1}s^{N-1}\partial^{\alpha}u(s
x+(1-s)y)ds\right|\leq ||\partial^{\alpha}u||_{L_{\infty}}
\end{equation*}
we obtain
    \begin{equation*}
        |\cR_{h,N}(x)|
\leq (h \sqrt{\cD})^{N} \sum_{|\beta|=N} \big\Vert \de^{\beta}
       u\big\Vert_{L^{\infty}}  \sum_{0< \alpha\leq \beta}
        \frac{|a_{\beta-\alpha}|}{\alpha!} \Vert \rho_{\alpha}(\cdot,\cD,\cH)
        \Vert_{L_{\infty}}
\end{equation*}
where
   \begin{equation*}
       \rho_{\alpha}(\xi,\cD,\cH)= \cD^{-n/2}\sum_{m\in\Z^{n}}\Big|
\Big(\frac{\xi-m}{\sqrt{\cD}}\Big)^{\alpha}
\cH \Big(\frac{\xi-m}{\sqrt{\cD}}\Big)\Big|.
   \end{equation*}
   In view of the decay condition there exist constants $c_{\alpha}$
   such that, for $\cD$ sufficiently large (\cite{MS4})
   \[\Vert \rho_{\alpha}(\cdot,\cD,\cH)
        \Vert_{L_{\infty}}\leq c_{\alpha}, \quad 0\leq |\alpha|\leq N.\]
Hence we obtain
     \begin{equation}\label{3.13}
        |\cR_{h,N}(x)|
\leq c\, (h \sqrt{\cD})^{N} \sum_{|\beta|=N} \big\Vert \de^{\beta}
       u\big\Vert_{L^{\infty}}
    \end{equation}
with a constant $c$ not depending on $h$, $\cD$ and $u$.
By \eqref{3.6}, this leads together with
(\ref{3.11}) and  (\ref{3.12})
to the estimate \eqref{3.14}.

\section{Construction of generating functions for arbitrary $N$}\label{S:QUA}

Here we describe a simple method for the construction of generating functions
$\cH$
satisfying the conditions (\ref{3.10}) for arbitrary
given $a_{\gamma}$ with $a_{0}=1$. Let us denote by $\cA=\{\cA_{\alpha\beta}\}$
the triangular matrix with the elements
    \[
\cA _{\alpha\beta}=
\begin{cases}
    a_{\beta-\alpha}\quad\> &\alpha\leq \beta,\cr
0 \quad &{\rm otherwise}
\end{cases}
\quad |\beta|,|\alpha|=0,\ldots,N-1.
\]
The dimension of $\cA$ is  $(N+n-1)!/((N-1)! \,n!)$.
Since $\det\cA=1$ there exists the inverse matrix
$\cA^{-1}=\{\cA_{\alpha\beta}^{(-1)}\}$ and
(\ref{3.10}) leads to the following conditions
\begin{equation}\label{3.10bis}
           \int_{\R^{n}} x^{\alpha} \cH(x)\, dx= \alpha ! \, \cA_{0
    \alpha}^{(-1)},\quad |\alpha|=0,\ldots,N-1.
           \end{equation}
These conditions can be rewritten as
\[
\partial^{\alpha} ( \cF \cH)(0)= \alpha !(-2\pi i)^{|\alpha|} \cA^{(-1)}_{0\alpha} ,\quad
|\alpha|=0,\ldots,N-1.
\]
Let us assume (see \cite{MS3})
\[
\cH(x)=\cP_{N}\Big(\frac{1}{2 \pi i} \frac{ \partial}{\partial x}\Big)\, \eta(x),
\]
where $\cP_{N}(t)$ is a polynomial of degree less than or equal to
$N-1$ and $\eta$ is a smooth function rapidly
decaying as $|x|\to \infty$ with $\cF\eta(0) \neq 0$.
Conditions (\ref{3.10bis}) give
\[
\partial^{\alpha}\big(\cP_{N}(\lambda) \cF \eta(\lambda)\big)(0)=(-2 \pi
i)^{|\alpha|}\alpha!\,\cA^{(-1)}_{0\alpha},\,
\quad |\alpha|=0,\ldots,N-1.
\]
We choose
$\cP_{N}$ as the
Taylor polynomial of order $(N-1)$
of the function
\[
Q(\lambda)=\sum_{|\alpha|=0}^{N-1} \frac{\cA_{0\alpha}^{(-1)}(-2\pi
i\lambda)^{\alpha}}{\cF\eta(\lambda)},\quad \lambda \in \R^{n}.
\]
Since
\[
\partial^{\beta}Q(0)=\sum_{\alpha\leq \beta}\cA_{0 \alpha}^{(-1)}(-2\pi i)^{|\alpha|}
\frac{\beta!}{(\beta-\alpha)!} \partial^{\beta-\alpha} (\cF\eta)^{-1}(0),
\]
where we use the notation
\[
\partial^{\beta-\alpha}(\cF\eta)^{-1}(0)=\partial^{\beta-\alpha}\Big(\frac{1}{\cF\eta(\lambda)}\Big)(0),
\]
we obtain
\[
\cP_{N}(\lambda)=\sum_{|\beta|=0}^{N-1}\frac{\lambda^{\beta}}{\beta!}
\sum_{\alpha\leq \beta}\cA_{0 \alpha}^{(-1)}(-2\pi i)^{|\alpha|}
\frac{\beta!}{(\beta-\alpha)!} \partial^{\beta-\alpha} (\cF\eta)^{-1}(0) \, .
\]
Therefore the equations
\begin{align*}
\partial^{\alpha} \big(\cP_{N}(\lambda) \cF \eta(\lambda)\big)(0)
&= \, \partial^{\alpha}(Q(\lambda)\cF\eta(\lambda))(0)\\
&= \, \partial^{\alpha} \Big(
\sum_{|\gamma|=0}^{N-1} {\cA_{0\gamma}^{(-1)}(-2\pi
i\lambda)^{\gamma}} \Big)(0)=(-2 \pi i)^{|\alpha|}\alpha! \, \cA^{(-1)}_{0\alpha}
\end{align*}
are valid for all $\alpha:0\leq|\alpha|\leq N-1$.
We have thus proved the following
\begin{theorem}\label{4.1T}
           Suppose that $\eta\in C^{N-1}(\R^{n})$  satisfies
\begin{align*}
       &|\eta(x)|\leq A\, (1+|x|)^{-K},\quad x \in \R^{n}, \quad K>N+n, \\
           &\itg_{\R^{n}} |x|^{N-1}|\partial^{\alpha}\eta(x)|dx<\infty,\quad 0 \leq
           |\alpha|\leq N-1 ,
\end{align*}
and $\cF\eta(0)\neq 0$.
Then the function
$$\cH(x)=\sum_{|\beta|=0}^{N-1}\sum_{\alpha \leq \beta}
\cA_{0 \alpha}^{(-1)}{ (-1)^{|\alpha|}}
\frac{\partial^{\beta-\alpha} (\cF\eta)^{-1}(0) }{ (\beta-\alpha)!(2 \pi
i)^{|\beta-\alpha|}}  \partial^{\beta}\eta(x)
$$
satisfies the conditions \eqref{3.10}.
\end{theorem}

Suppose that $\eta(x)$ is radial, that is $\eta(x)=\psi(r), r=|x|$. Then
$\partial ^{\alpha}\cF \eta(0)=0$ for any $\alpha=(\alpha_{1},\ldots,\alpha_{n})$
containing at least one odd $\alpha_{i}$ and
we obtain the formula
\[
\cH(x)=\sum_{|\beta|=0}^{N-1}(-\partial)^{\beta}\eta(x)\,\sum_{2\gamma \leq \beta}
\cA_{0, \beta-2\gamma}^{(-1)}
\frac{\partial^{2\gamma} (\cF\eta)^{-1}(0) }{ (2\gamma)!(-4 \pi^{2})^{|\gamma|}}
\, .
\]
Let us consider the
special case  $\eta(x)=\pi^{-n/2} {\ee}^{-|x|^{2}}$
with $\cF\eta(\lambda)={\ee}^{-\pi^{2}|\lambda|^{2}}$.
Denoting by $H_{\beta}$ the Hermite polynomial of $n$ variables defined by
\[
H_{\beta}(t)={\ee}^{|t|^{2}}(-{\partial}
)^{\beta}{\ee}^{-|t|^{2}}
\]
we derive $\partial^{\gamma}(\cF \eta)^{-1}(0)=(-\pi i)^{|\gamma|}H_{\gamma}(0)$.
Note that
$H_{\gamma}(0)=0$  if $\gamma$ has  odd components,
otherwise $H_{2 \gamma}(0)=(-1)^{\gamma}{(2\gamma)!}/{\gamma!}$.
Hence
\begin{equation}\label{4.1}
           \cH(x)=\pi^{-n/2}\sum_{|\beta|=0}^{N-1}H_{\beta}(x){\ee}^{-|x|^{2}}
           \sum_{2\gamma\leq \beta}\frac{(-1)^{|\gamma|}}{\gamma!
4^{|\gamma|}}
{\cA_{0\, ,\beta-2\gamma}^{(-1)}}.
\end{equation}
Assuming $a_{\gamma}=\delta_{|\gamma| 0}$ and $N=2M$ we find out
\begin{align*}
\cH(x)= & \, \pi^{-n/2}\sum_{j=0}^{M-1}\frac{(-1)^{j}}{4^{j}}
\sum_{|\beta|=j}\frac{H_{2\beta}(x)}{\beta!}
\, \ee^{-|x|^{2}}\\
=& \, \pi^{-n/2}\sum_{j=0}^{M-1}\frac{(-1)^{j}}{4^{j}j!}\Delta^{j}
{\ee}^{-|x|^{2}}=\pi^{-n/2}L_{M-1}^{(n/2)}(|x|^{2}){{\ee}^{-|x|^{2}}}
\,,
\end{align*}
where $L_{k}^{(\gamma)}$ are
the generalized Laguerre polynomials,
which
are
defined by
       \begin{equation*}
L_{k}^{(\gamma)}(y)=\frac{\ee^{\,y} y^{-\gamma}}{k!} \, \Big(
{\frac{d}{dy}}\Big)^{k} \!
\left(\ee^{\,-y} y^{k+\gamma}\right), \quad \gamma > -{\rm 1} \, .
       \end{equation*}
In this case
we obtain the classical generating function
$\eta(x)=L^{(n/2)}_{M-1}(|x|^{2}){\ee}^{-|x|^{2}}$ (see
\cite{MS3,Sc}).

\section{Examples}\label{example}
\begin{example}
If $N=2$, then formula \eqref{4.1}  gives
\begin{equation}
       \cH(x)=\pi^{-n/2}\Big(1-2 \sum_{|\beta|=1}
       a_{\beta}x^{\beta}\Big){\ee}^{-|x|^{2}}.
\end{equation}
For the one-dimensional case the corresponding quasi-interpolant is
\[
M_{a} u(x)=(\pi \cD)^{-1/2}\sum_{m\in \Z} (u(hm)-h \sqrt{\cD}\, a
\,u'(h m))(1-2
a \,\frac{x-h m}{h \sqrt{\cD}}\,){\ee}^{-(x-hm)^{2}/(h^{2}\cD)}\, .
\]
For $u(x)=\cos x$, $x \in \R$, the difference between $u(x)$ and
$M_{a} u(x)$ is plotted in \ref{fig8}
by taking $h=0.1$, $\cD=2$ for different values of $a$.

\begin{figure}\begin{center}
\epsfig{file=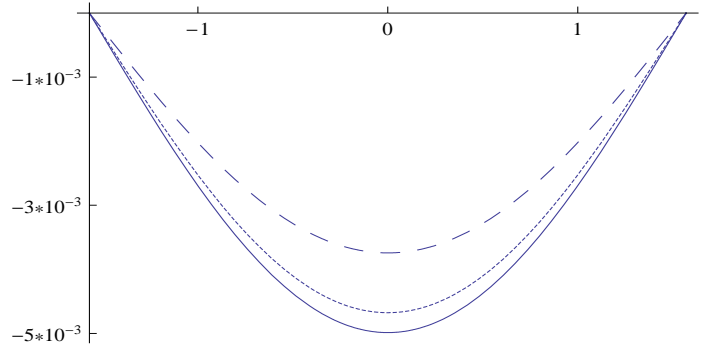,width=3.6in,height=1.8in}
\caption{The graphs of $(M_{a} -I)\cos x $ when $a=0$ (solid
line), $a=1/8$ (dotted line) and $a=1/4$ (dashed line).}
\label{fig8}
\end{center}
\end{figure}

\end{example}

\begin{example}
       Now we are looking for quasi-interpolants of order $\cO(h^{4})$.
In the one-dimensional case
\[
\cA^{(-1)}_{01}=-a_{1},\>
\cA_{02}^{(-1)}=a_{1}^{2}-a_{2},\>
\cA^{(-1)}_{03}=-a_{1}^{3}+2 a_{1}a_{2}-a_{3}
\]
and formula (\ref{4.1}) gives
\[\begin{split}
   \cH(x)=&
\pi^{-1/2}{\rm
e}^{-x^{2}}  [(3/2-\cA_{02}^{(-1)})+(5\cA_{01}^{-1}-12\cA_{03}^{(-1)})x \\
&+ (4\cA_{02}^{(-1)}-1)x^{2}+(8\cA_{03}^{(-1)}-2\cA_{01}^{(-1)})x^{3}].
\end{split}\]

If $a_{1}=a_{2}=a_{3}=0$, then we get the  classical generating function
$\eta(x)=\pi^{-1/2}{\rm
e}^{-x^{2}}(3/2-x^{2})$.
Assuming $a_{3}=a_{2}=0,\, a_{1}=\pm 1/2$ we
obtain
\[
\cH(x)=\pi^{-1/2}(1\pm x){\ee}^{-x^{2}}
\]
and the quasi-interpolant  of order $\cO(h^{4})$:
\[
M_{1}u(x)=(\pi \cD)^{-1/2}\sum_{m\in \Z}(u(h m)\pm \frac{h \sqrt{\cD}}{2}
u'(hm))(1\pm \frac{x-h m}{h \sqrt{\cD}})\, {\ee}^{-(x-h m)^2/(h^2 \cD)} \,.
\]
With the choice $a_{1}=a_{3}=0,\, a_{2}=-1/4$ we
obtain
$\cH(x)=\pi^{-1/2}{\ee}^{-x^{2}}$
and the quasi-interpolant  of order $\cO(h^{4})$:
\[
M_{2}u(x)=(\pi \cD)^{-1/2}\sum_{m\in \Z}(u(h m)- \frac{h^{2} {\cD}}{4}
u''(hm))\, {\ee}^{-(x-h m)^2/(h^2 \cD)} \, .
\]

In the Figures \ref{fig5}, \ref{fig6}, \ref{fig7} we show the error
graphs for the
approximation of $u=\cos(x)$ with the quasi-interpolants $Mu, M_{1}u,
M_{2}u$, respectively. We have taken $\cD=2$, $h=0.05$, $h=0.1$ and
$h=0.2$.  The
case of
smaller $h$ gives different pictures: it is clearly visible that the
error oscillates very fast
with $M u$ and $M_{1}u$.

\begin{figure}\begin{center}
    \subfigure[]{
\epsfig{file=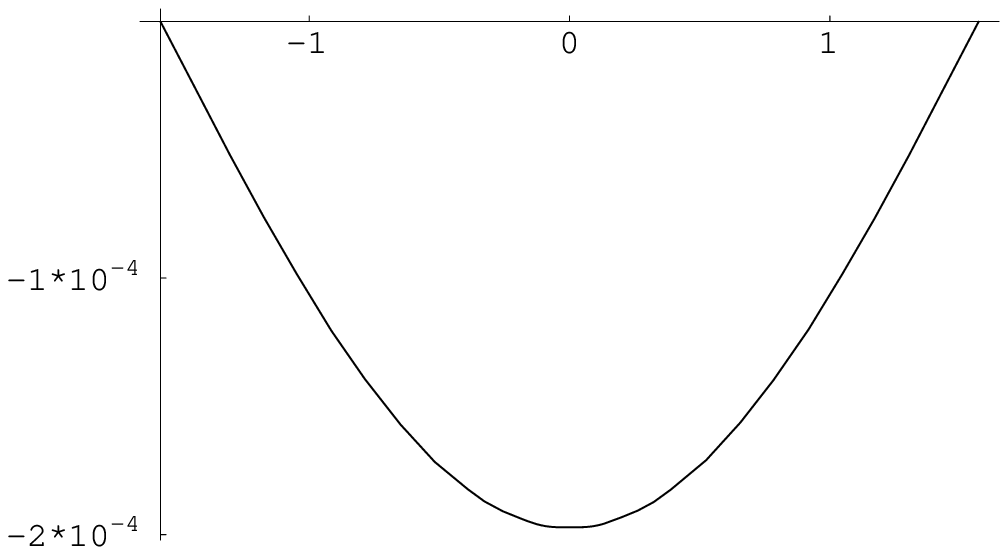,width=1.6in,height=1.8in}}
 \subfigure[]{
\epsfig{file=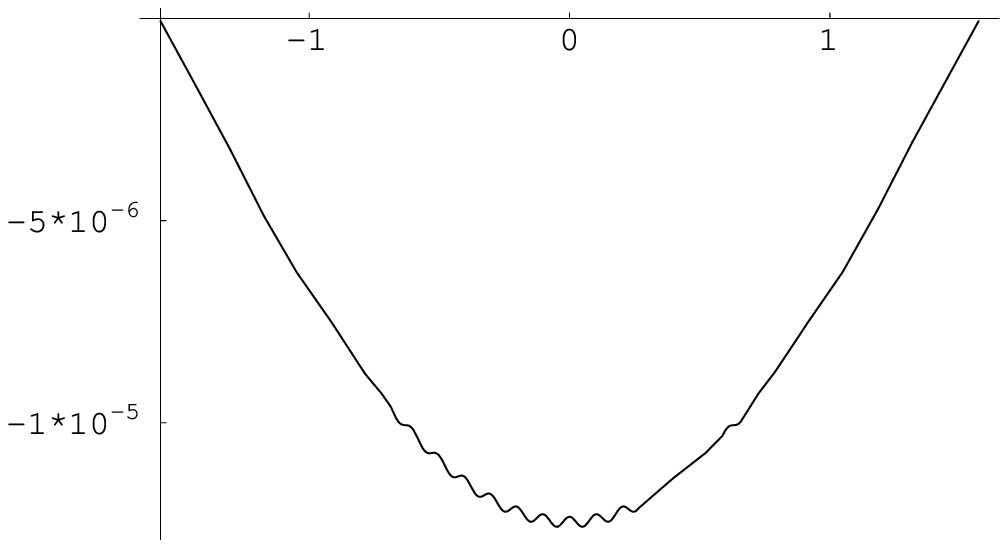,width=1.6in,height=1.8in}}
 \subfigure[]{
\epsfig{file=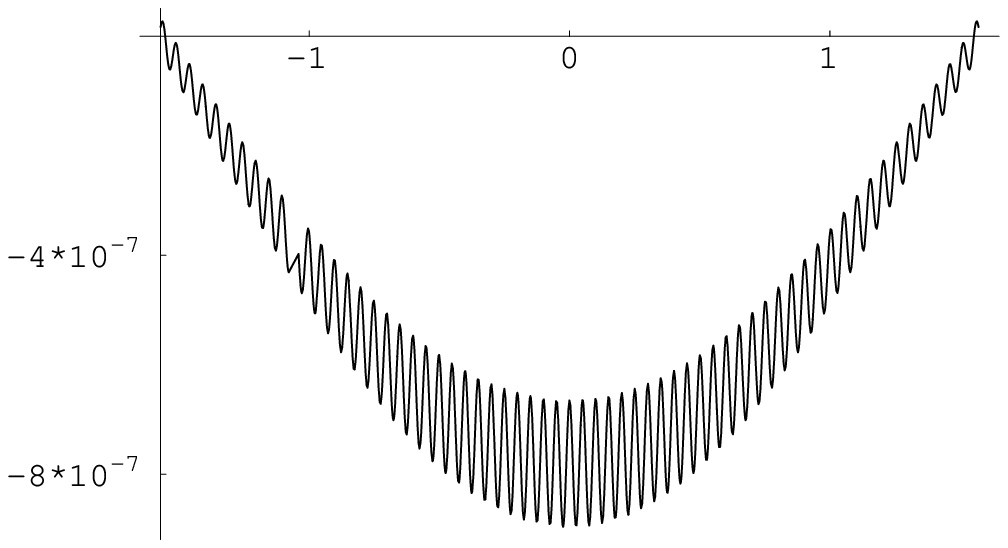,width=1.6in,height=1.8in}}
\caption{The graphs of $(M-I)\cos(x)$ by assuming $h=0.2$
(a), $h=0.1$ (b) and $h=0.05$ (c).}
\label{fig5}
\end{center}
\end{figure}

\begin{figure}\begin{center}
     \subfigure[]{
\epsfig{file=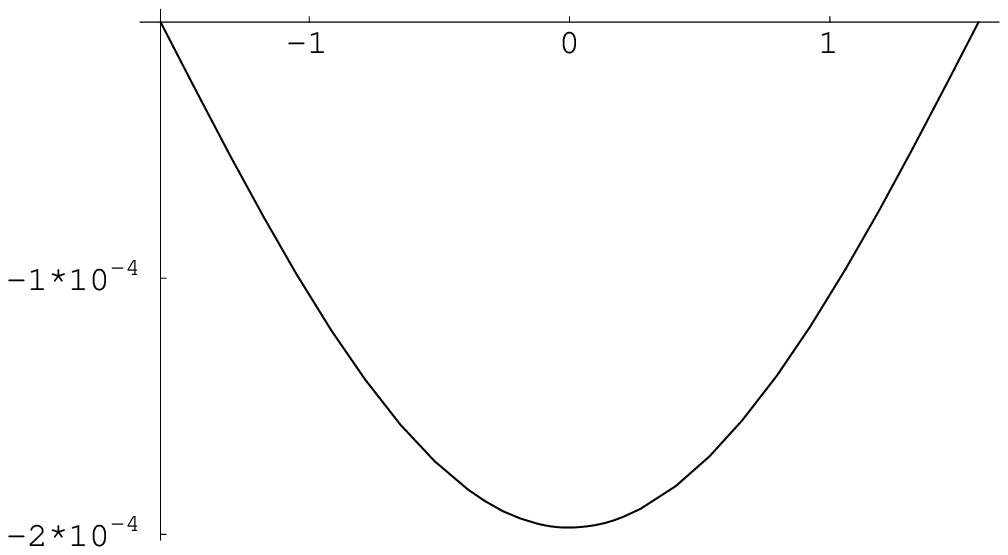,width=1.6in,height=1.8in}}
 \subfigure[]{
\epsfig{file=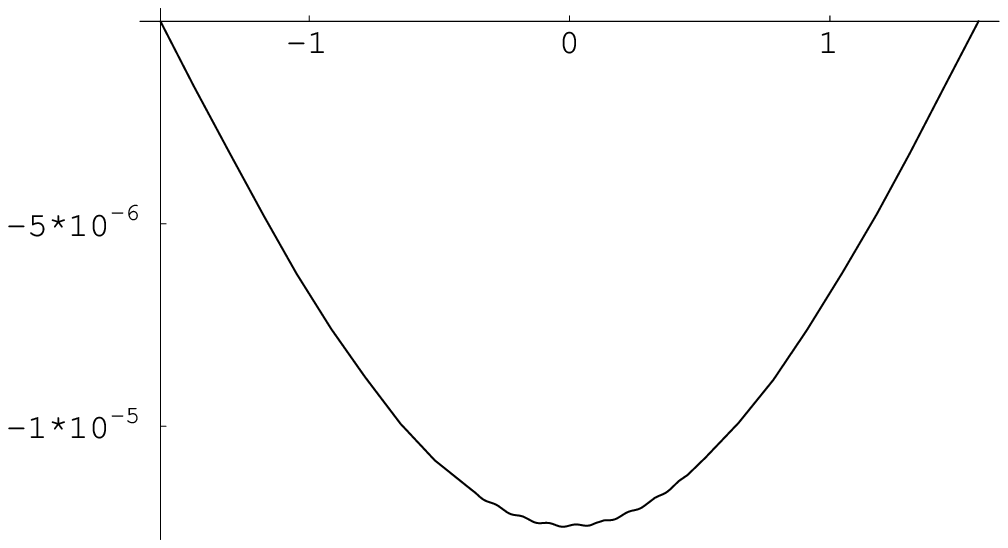,width=1.6in,height=1.8in}}
 \subfigure[]{
\epsfig{file=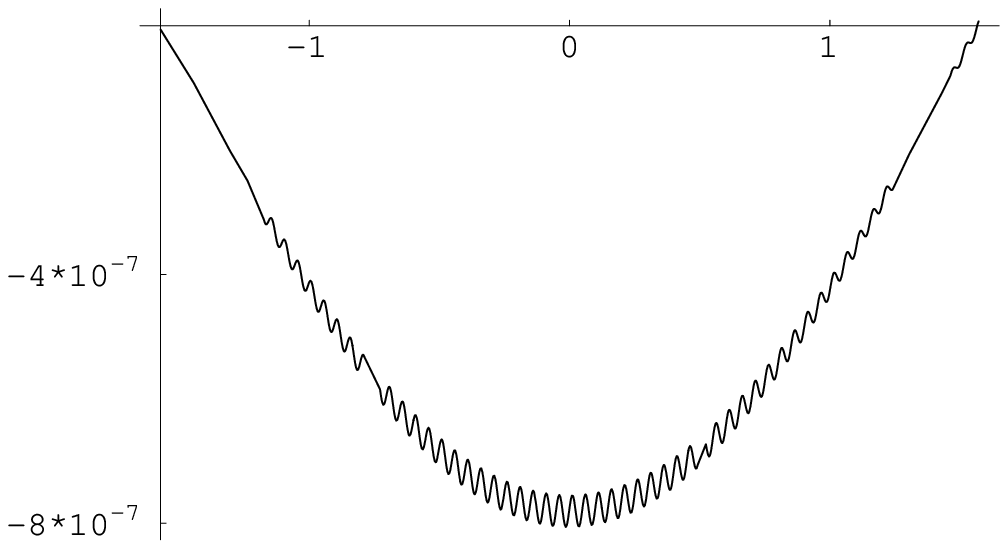,width=1.6in,height=1.8in}}
\caption{The graphs of $(M_{1}-I)\cos(x)$ by assuming $h=0.2$
(a), $h=0.1$ (b) and $h=0.05$ (c).}
\label{fig6}
\end{center}
\end{figure}

\begin{figure}\begin{center}
     \subfigure[]{
\epsfig{file=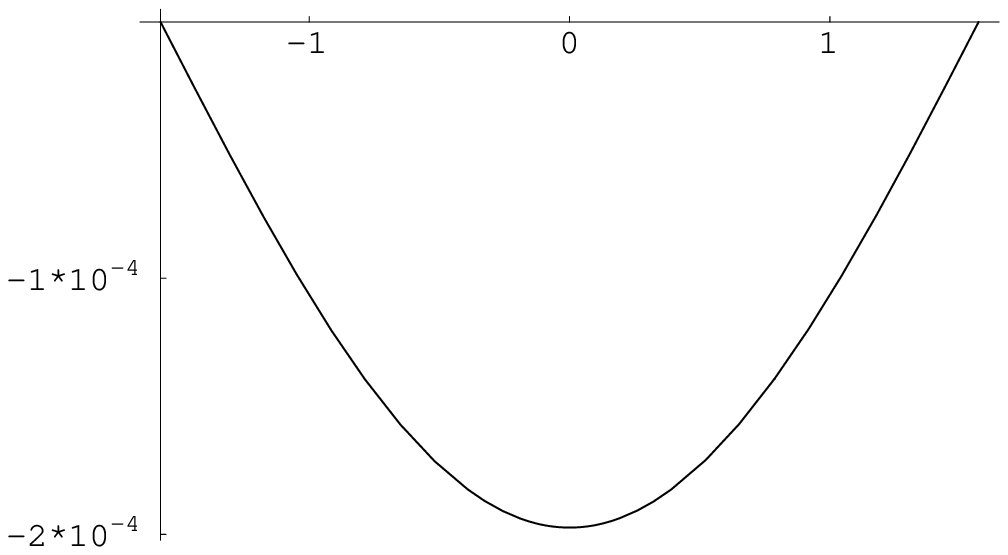,width=1.6in,height=1.8in}}
 \subfigure[]{
\epsfig{file=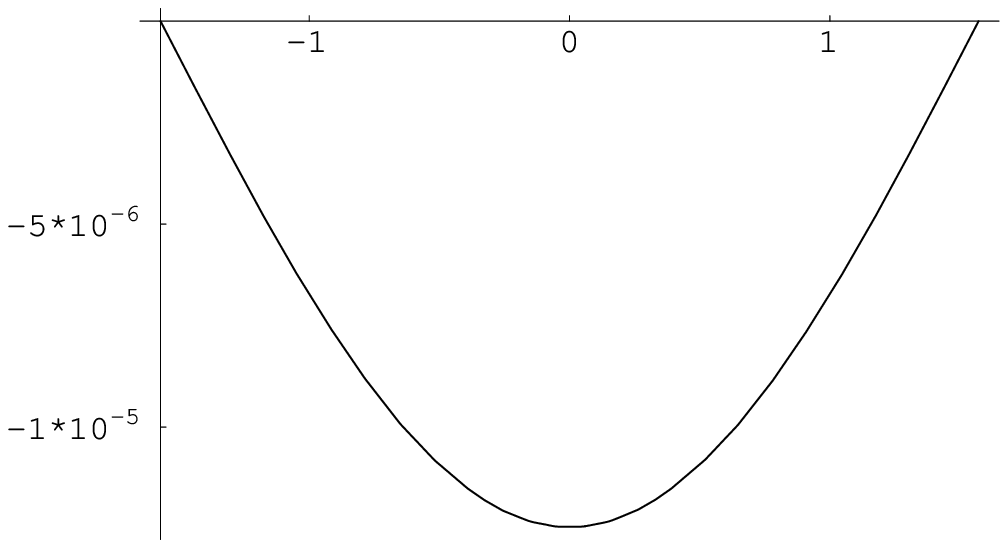,width=1.6in,height=1.8in}}
 \subfigure[]{
\epsfig{file=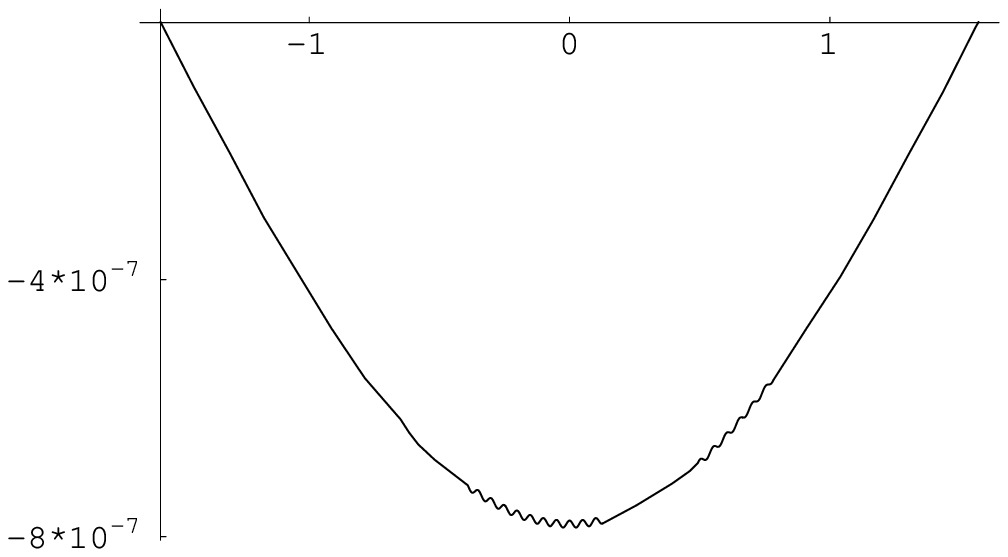,width=1.6in,height=1.8in}}
\caption{The graphs of $(M_{2}-I)\cos(x)$ by assuming $h=0.2$
(a), $h=0.1$ (b) and $h=0.05$ (c).}
\label{fig7}
\end{center}
\end{figure}

If $n>1$,
by taking $a_{2\beta}=a$ and $a_{\beta}=0$ otherwise, we
obtain
the radial generating function
\begin{align*}
  \cH_{a}(x)=& (1-\sum_{|\beta|=1}
(a_{2\beta}+1/4)(4 x^{2\beta}-2)) {\ee}^{-|x|^{2}}\pi^{-n/2} \\
=&(1+n(2 a+1/2)-(1+4a)|x|^{2} ){\ee}^{-|x|^{2}}{\pi^{-n/2}}
\end{align*} 
for the approximate approximation of order $\cO(h^{4})$
\[
M_{a}u(x)=(\pi \cD)^{-n/2}\sum_{m\in \Z^{n}}(u(h m)+h^{2}\, \cD\, a \,\Delta
u(hm))\cH_{a}\big( \frac{x-hm}{h\sqrt{\cD}}\big).
\]
\end{example}

\begin{example}\label{4.3}
Let $n>1$ and set in (\ref{3.1}) $N=2M$ and
\begin{align} \label{coefLap}
\begin{cases} \;
a_{\gamma}\,=0 \, ,\quad {\rm if}\; \gamma \; \mbox{has odd
components,}\\[3mm]
\;  \displaystyle{a_{2\gamma}=\frac{(-1)^{|\gamma|}}{\gamma! 4^{|\gamma|}}},\quad 0\leq
|\gamma|\leq M-1.
\end{cases}
\end{align}
Keeping in mind that
\begin{equation}\label{lapl}
\sum_{|\gamma|=s}\frac{s!}{\gamma!} \, \de^{2\gamma}u(x)=\Delta^{s}u(x)
\end{equation}
we obtain
\[
\cQ\left((-h
           \sqrt{\cD})\partial \right)u =
           \sum_{j=0}^{M-1}\frac{(-1)^{j}(h
           \sqrt{\cD})^{2j}}{4^{j}}\sum_{|\gamma|=j}\frac{\partial^{2
           \gamma}}{\gamma!}u=\sum_{j=0}^{M-1}\frac{(-1)^{j}(h^2 
\cD)^j}{j!4^{j}}\Delta^{j}u.
\]
The function $\cH(x)=\pi^{-n/2}{\ee}^{-|x|^{2}}$ satisfies conditions
(\ref{3.10}). In fact
\[
\frac{1}{\pi^{n/2}}\itg_{\R^{n}}x^{\alpha} \ee^{-|x|^2} dx=\frac{(-\pi i)^{\alpha}}{(2
\pi)^{\alpha}}H_{\alpha}(0)=
       \begin{cases}
\; 0, & \mbox{if } \alpha \mbox{ has odd components,} \\[2mm]
\ds \frac{(2 \gamma)!}{\gamma!2^{2|\gamma|}} ,
& \mbox{if } \alpha=2\gamma \, .
\end{cases}
\]
Then the equations (\ref{3.10}) are valid for all $\beta: 0 \leq |\beta|\leq M-1$
because of
\[
\sum_{\alpha\leq \beta}\frac{a_{2(\beta-\alpha)}}{\alpha!2^{2|\alpha|}}=
\sum_{\alpha\leq \beta} \frac{(-1)^{|\beta-\alpha|}}{(\beta-\alpha)!4^{|\beta-\alpha|}}
\frac{1}{\alpha!4^{|\alpha|}}=\delta_{|\beta| 0},\quad |\beta|=0,\ldots,M-1.
\]
Therefore, a  general approximation of order $N=2M$ is given by
\begin{equation}\label{hope}
       M^{(N)}u(x)=(\pi \cD)^{-n/2}
\sum_{m\in\Z^{n}}\sum_{s=0}^{M-1} (h \sqrt{\cD})^{2s}
\frac{(-1)^{s}}{s!4^{s}}\Delta^{s}u(hm)\, {\ee}^{-|x-h m|^2/(h^2 \cD)} \, .
\end{equation}
If $M=2$ then we find the ``fourth order formula''
\begin{equation*}
M^{(4)}u(x)=(\pi \cD)^{-n/2}\sum_{m\in \Z^{n}} \big[u(h m)- \frac{h^2
         \cD}{4}\Delta u(h m)
          \big]   {\ee}^{-|x-h m|^2/(h^2 \cD)} \, .
\end{equation*}
For $M=3$ we obtain the ``sixth order formula''
\[
M^{(6)}u(x)=(\pi \cD)^{-n/2}\!\!\sum_{m\in \Z^{n}}\!\!\! \big[u(h m)-
        \frac{h^2\cD}{4}\Delta u(h m)
        +\frac{h^4 \cD^{2}}{32}\Delta^{2}u(hm)
          \big]  {\ee}^{-|x-h m|^2/(h^2 \cD)} .
\]
Note the additive
structure of the formula (\ref{hope})
\[
M^{(N+2)}u(x)=M^{(N)}u(x)+ \frac{(h
          \sqrt{\cD})^{N}(-1)^{M}}{(\pi \cD)^{n/2}M! 4^{M}}\!\!\!\!\sum_{m\in \Z^{n}}
          \Delta^{M} u(h m) \, {\ee}^{-|x-h m|^2/(h^2 \cD)}.
\]
\end{example}

\begin{example}\label{generalcase}
 We consider the second order partial differential operator
\[
\cB u=\sum_{i,k=1}^{n} b_{ik} \partial_{i}\partial_{k} u \, ,
\]
where the matrix $B=\{b_{ik}\}\in\R^{n\times n} $ is symmetric and
positive definite.
Define the quasi-interpolant
\begin{equation}\label{gen}
       M u(x)=\frac{(\det B)^{-1/2}}{(\pi \cD)^{n/2}}
\sum_{m\in\Z^{n}}\sum_{s=0}^{M-1}
\frac{(-h^2 \cD )^{s}}{s!4^{s}} \, \cB ^{s}u(hm)\,\ee^{\, -\langle
   B^{-1}(x-hm),x-hm\rangle/(h^2\cD)} .
\end{equation}
Assume $C$ such that
      $B^{-1}=C^{T}C$.
If we consider the linear transformation $\xi=C\,x$ and introduce
$U(\xi)=u(x)$, in the
new coordinates we have (\cite[p.42]{Pe})
\[
\cB u(x)=\Delta U(\xi).
\]
Then \eqref{gen} will take the form
\begin{align*}
       M u(x)={\cM} U(\xi)=
\frac{\det C}{(\pi \cD)^{n/2}}
\sum_{m\in\Z^{n}}\sum_{s=0}^{M-1}
\frac{(-h^2 \cD)^{s}}{s!4^{s}}\Delta^{s} U(h\,C\,m)\,
\ee^{\, -|\xi-h\,C\,m|^2/(h^2\cD)}.
\end{align*}
Keeping in mind \eqref{lapl} and \eqref{3.2} we rewrite
\[
{\cM}U(\xi)=\sum_{|\beta|=0}^{2M-1}(-h
\sqrt{\cD})^{|\beta|}\de^{\beta}U(\xi) \,
\cE_{\beta}\big(\frac{\xi}{h},\cD\big)+R_{h}(\xi) \,
\]
where
\begin{align*}
\cE_{\beta}(\xi,\cD) &=\sum_{2\gamma\leq\beta}\frac{(-1)^{|\gamma|}}{\gamma!(\beta-2\gamma)!4^{|\gamma|}}
\Sigma_{\beta-2\gamma}({\xi},\cD)
\intertext{with}
\Sigma_{\alpha}(\xi,\cD) &=\frac{\det C}{(\pi
\cD)^{n/2}}\sum_{m\in\Z^{n}}\left(\frac{\xi-C\,m}{\sqrt{\cD}}\right)^{\alpha}
\!\!{\ee}^{-|\frac{\xi-C\,m}{ \sqrt{\cD}}|^{2}} \, ,
\end{align*}
and
\begin{align*}
&R_{h}(\xi)= \, \frac{(h\sqrt{\cD})^{2M}\det C}{(\pi
\cD)^{n/2}} \\
&\times \sum_{|\beta|=2M}
\sum_{2\gamma\leq\beta}\frac{(-1)^{|\gamma|}}{4^{|\gamma|}\gamma!(\beta-2\gamma)!}
\sum_{m\in\Z^{n}}U_{\beta}(\xi,hCm)
\left(\frac{\xi-hCm}{h\sqrt{\cD}}\right)^{\beta-2\gamma}
{\ee}^{-|\frac{\xi-hCm}{ h\sqrt{\cD}}|^{2}}.
\end{align*}
Now we use Poisson's summation formula on affine grids (see
\cite[p.23]{MSbook})
\begin{equation}\label{poisson}
\begin{split}
   \frac{\det C}{
\cD^{n/2}}\sum_{m\in\Z^n}\hskip-3pt
\Big(\frac{\xi-C\,m}{\sqrt{\cD}}\Big)^{\delta}
\eta\Big(\frac{\xi-C\,m}{\sqrt{\cD}}\Big)&=   \\
\Big(\frac{i}{2\pi}\Big)^{|\delta|}&\sum_{\nu\in\Z^{n}}\de^{\delta}\cF\eta
(\sqrt{\cD}C^{-T}\nu)\ee^{2\pi i
\langle  \xi,C^{-T}\nu\rangle },
\end{split}
\end{equation}
where we denote by $C^{-T}= (C^T)^{-1}$.
In our case $\eta(x)=\pi^{-n/2} \ee^{-|x|^2}$
and
$\cF\eta(\lambda)=\ee^{-\pi^{2}|\lambda|^{2}}$. Since
$\de^{\delta}\cF\eta(\lambda)=(-\pi)^{\delta}H_{\delta}(\pi\,\lambda) {\rm
e}^{-\pi^{2}|\lambda|^{2}} $
we obtain
\[
\de^{\delta}\cF\eta(0)=\begin{cases}
0\, , & \mbox{if } \delta \mbox{ has odd components,} \\[3mm]
\ds \pi^{2\gamma}\frac{(-1)^{|\gamma|}(2\gamma)!}{\gamma!}\, , & \mbox{if } \delta=2\gamma \, .
  \end{cases}
\]
Formula \eqref{poisson} applied to $\Sigma_{\alpha}$ gives
\[
\Sigma_{\alpha}(\xi,\cD)=
\begin{cases}
\ds\Big(\frac{i}{2\pi}\Big)^{|\alpha|}\sum_{\nu\neq 0}\de^{\alpha}\cF\eta
(\sqrt{\cD}C^{-T}\nu)\ee^{2\pi i
\langle\xi,C^{-T}\nu\rangle } , \, \alpha \mbox{ has odd 
components,}
\\[4mm]
\ds
\frac{\alpha!}{2^{|\alpha|} \gamma!}+
\Big(\frac{i}{2\pi}\Big)^{|\alpha|}\ds\sum_{\nu\neq 0}\de^{\alpha}\cF\eta
(\sqrt{\cD}C^{-T}\nu)\ee^{2\pi i
\langle\xi,C^{-T}\nu\rangle }, \, \alpha=2\gamma .
  \end{cases}
\]
We deduce that
\[
\begin{split}
     \cE_{\beta}(\xi,\cD)&=\delta_{|\beta|0} \\
&+\sum_{2\gamma\leq\beta}
\frac{(-1)^{|\gamma|}}{\gamma!(\beta-2\gamma)!\, 4^{|\gamma|}}
\Big(\frac{i}{2\pi}\Big)^{|\beta-2\gamma|}\sum_{\nu\neq 0}\de^{\beta-2\gamma}\cF\eta
(\sqrt{\cD}C^{-T}\nu){\rm e}^{2\pi i
\langle\xi,C^{-T}\nu\rangle }
\end{split}
\]
and
\[
|\cE_{\beta}(\xi,\cD)-\delta_{|\beta|0}|\leq
\sum_{2\gamma\leq\beta}\frac{\pi^{|2\gamma-\beta|}}{\gamma!(\beta-2\gamma)!\, 2^{|\beta|}}
\sum_{\nu\neq 0} \big|\de^{\beta-2\gamma}\cF\eta
(\sqrt{\cD}C^{-T}\nu)\big|.
\]
By repeating the same arguments used in the proof of Theorem \ref{UNO}
we derive that  $|\cE_{\beta}(\xi,\cD)-\delta_{|\beta|0}|<\e_{1}$ for prescribed
$\e_{1}>0$ and sufficiently large $\cD$, and the remainder is bounded by
\[
|R_{h}(\xi)|\leq c_{B}
(h\sqrt{\cD})^{2M}\sum_{|\beta|=2M}\|\de^{\beta}U\|_{L^{\infty}}
\]
with the constant $c_{B}$ independent of $U,h,\cD$.
Hence
\begin{align*}
|M u(x) -u(x)|&=|{\cM} U(\xi)-U(\xi)|\\
&  \leq \e_{1}\sum_{|\beta|=0}^{2M-1}(h
\sqrt{\cD})^{|\beta|}|\de^{\beta}U(\xi)|
+c_{B}
(h\sqrt{\cD})^{2M}\sum_{|\beta|=2M}||\de^{\beta}U||_{L^{\infty}}\\
&  \leq
\e\sum_{|\beta|=0}^{2M-1}(h
\sqrt{\cD})^{|\beta|}|\de^{\beta}u(x)|
+C_{B}
(h\sqrt{\cD})^{2M}\sum_{|\beta|=2M}||\de^{\beta}u||_{L^{\infty}} \, ,
\end{align*}
where $C_{B}$ depends only on the matrix $B$.
Therefore the quasi-interpolant
\eqref{gen} approximates $u$ with the order $\cO((h\sqrt{\cD})^{2M})$
up to the saturation error.

\end{example}

\section{An application of formula (\ref{hope})}\label{appl}
Here we consider the approximation of harmonic functions.
Suppose that $\Delta u=0$  in some domain $\Omega \subseteq \R^n$. 
Then for any $N=2M$ and $x \in \Omega$
the Hermite  quasi-interpolant \eqref{hope}  has  the simple form
\begin{align} \label{harmonic}
       M^{(N)}u(x)=M u(x)=(\pi \cD)^{-n/2}
\sum_{hm\in \Omega} u(hm)\, {\ee}^{-|x-h m|^2/(h^2 \cD)} \, ,
\end{align}
i.e.,  it coincides with the well known quasi-interpolation formula
of second order. However, Theorem \ref{UNO} indicates higher 
approximation rates.
This will be studied here in more detail.

First we consider the case  $\Omega = \R^{n}$.
Then
\[
u(\xi)=\sum_{|\beta|=0}^{\infty}\frac{\de^{\beta}u(x)}{\beta!}\,(\xi-x)^{\beta},\quad
  \xi\in\R^{n} 
\]
and  the series  converges absolutely in $\R^{n}$.
Moreover, $u$ has the
analytic extension
\begin{align} \label{anext}
\tilde u (\zeta)=\sum_{|\beta|=0}^{\infty}\frac{\de^{\beta}u(x)}{\beta!}\, (\zeta-x)^{\beta},\quad
     \zeta\in\C^{n} \, ,
\end{align}
cf. e.g.  \cite{Se, Fi}.
Using formula \eqref{3.6} for the quasi-interpolant with the
generating function \eqref{hope}
we obtain
\begin{align} \label{4.6}
M u(x)-u(x) = \cE_{h,2M}(x)+\cR_{h,2M}(x) \, ,
\end{align}
where
\begin{align*}
&\cE_{h,2M}(x) =
u(x) \Big((\pi \cD)^{-n/2}  \sum_{m\in \Z^n}
{\ee}^{-|x-hm|^2/(h^2\cD)} -1\Big) \\
&\quad +(\pi \cD)^{-n/2} \sum_{|\beta|=1}^{2M-1} (-h \sqrt{\cD})^{|\beta|} 
\de^{\beta}u(x)
   \sum_{m\in \Z^n} {\ee}^{-|x-hm|^2/(h^2\cD)}
\sum_{\alpha\leq \beta} \frac{a_{\beta-\alpha}}{\alpha!}
\Big( \frac{x- hm}{h\sqrt{\cD}}\Big)^{\alpha}
\end{align*}
constitutes the saturation error and
the remainder term has the form
\begin{equation} \label{rem}
\begin{split}
&\cR_{h,2M}(x)=
(-h \sqrt{\cD})^{N}
  (\pi \cD)^{-n/2} \\
&\quad \times \sum_{|\beta|=2M}
   \sum_{m\in \Z^n} U_{\beta}(x,h m) {\ee}^{-|x-hm|^2/(h^2 \cD)}
   \sum_{0<\alpha\leq\beta} \frac{a_{\beta-\alpha}}{\alpha!}
\Big( \frac{x-h m}{h \sqrt{\cD}}\Big)^{\alpha}\,.
\end{split}
\end{equation}
From \eqref{coefLap} we see that
\begin{align*}
   \sum_{\alpha\leq\beta} \frac{a_{\beta-\alpha}}{\alpha!}
\,x^{\alpha}
=  \sum_{2\gamma\leq\beta}\frac{a_{2 \gamma}}{(\beta-2\gamma)!}
\,x^{\beta-2\gamma}
= \sum_{2\gamma\leq\beta}\frac{(-1)^{|\gamma|}}{\gamma! \, (\beta-2\gamma)! \, 2^{2|\gamma|}}
\, x^{\beta-2\gamma} \, ,
\end{align*}
which by using the representation of Hermite polynomials
\begin{align*}
H_k(\tau)=
\sum_{0\le 2j \le k} \frac{(-1)^j k!}{j! \, (k-2j)!}\,  (2\tau)^{k-2j} \,,
\end{align*}
shows that
\begin{align*}
   \sum_{\alpha\leq\beta} \frac{a_{\beta-\alpha}}{\alpha!}
\, x^{\alpha}
= \frac{1}{\beta! \, 2^{|\beta|}} H_\beta (x) \,.
\end{align*}
Hence we obtain
\begin{align*}
\cE_{h,2M}(x) =
\sum_{|\beta|=0}^{2M-1}
\Big( \hskip-3pt- \frac{h \sqrt{\cD}}{ 2}\Big)^{|\beta|}
   \frac{\de^{\beta}u(x)}{\beta! }\, \sigma_\beta \Big(\frac{x}{h}, \cD\Big)
\end{align*}
with the functions
\begin{equation} \label{defsigma}
\begin{split}
\sigma_0 (x, \cD) &=(\pi \cD)^{-n/2}  \sum_{m\in \Z^{n}}
{\ee}^{-|x-m|^2/\cD} -1 \, ,\\
\sigma_\beta (x, \cD) &= (\pi \cD)^{-n/2} \sum_{m \in \Z^n} H_\beta
\Big( \frac{x-m}{ \sqrt{\cD}}\Big)\, {\ee}^{-|x-m|^2/\cD}
\,,\; |\beta|=1,\ldots,2M-1 \, .
\end{split}
\end{equation}
It follows from the definition of $H_\beta$ and from Poisson's summation
formula that
\[
\sigma_\beta (x, \cD)= (-1)^{|\beta|} 
\cD^{|\beta|/2}\de^{\beta}\sigma_0 (x, \cD)=
(-2 \pi i )^{|\beta|} \cD^{|\beta|/2}
\!\!\!\sum_{m \in \Z^{n}\setminus\{0\}}\!\!\!\! m^\beta \ee^{ \,- \pi^2 \cD|m|^2}
\ee^{ \,2 \pi i \langle
    m,x\rangle}.
\]
Thus
the saturation error can be expressed as
\begin{align} \label{saturN}
\begin{split}
\cE_{h,2M}(x)
&= \sum_{|\beta|=0}^{2M-1}
   \frac{\de^{\beta}u(x)}{\beta! }\Big(\frac{h \cD}{2} \Big)^\beta
  \, \de^{\beta}\,\sigma_0 \Big(\frac{x}{h}, \cD\Big)\\
&=
\sum_{m \in \Z^{n}\setminus\{0\}}  \tilde u_{2M}(x+i \pi h \cD m)
\ee^{ \,- \pi^2 \cD|m|^2} \ee^{
\,2 \pi i \langle
    m,x\rangle/h}
  \, ,
\end{split}
\end{align}
where
\[
\tilde u_N(\zeta)=\sum_{|\beta|=0}^{N-1}\frac{\de^{\beta}u(x)}{\beta!}\, (\zeta-x)^{\beta}
\]
is the Taylor polynomial of the analytic extension $\tilde u$.
Note that \eqref{saturN} is valid for any $M$.
\begin{theorem}\label{thmarm}
Suppose that the harmonic in $\R^n$ function $u$ is such that
the series
\begin{align} \label{asssqu}
\sum_{|\beta|=0}^{\infty}\frac{\de^{\beta}u(x)}{\sqrt{\beta!}}\, y^{\beta}
\end{align}
converges absolutely for any $y \in \R^n$.
If $\sqrt{\cD}h <1$, then the quasi-interpolant \eqref{harmonic}
approximates $u$ with
\begin{align*}
M u(x)-u(x) = \lim_{M\to \infty}\cE_{h,2M}(x)=
\sum_{m\in \Z^{n}\setminus\{0\}} \tilde u(x+ \pi i h \cD m)
   \ee^{ \,- \pi^2 \cD|m|^2}  \ee^{ \, 2 \pi i \langle
    m,x\rangle/h} ,
\end{align*}
where $\tilde u$ is the analytic extension of $u$ onto $\C^n$.
\end{theorem}
{\bf Proof.}
We have to show that $|\cR_{h,2M}(x)| \to 0$ as $M \to \infty$.
To estimate
\eqref{rem} we rewrite
\begin{align*}
\sum_{0\le 2j < k} \frac{(-1)^j }{j! \, (k-2j)! \, 2^{2j}}\,  \tau^{k-2j}
&= \frac{1}{2^k \, k!  } \big( H_k(\tau) - H_k(0) \big)\\
&= \frac{1}{2^{k-1} \, (k-1)!  } \itg_0^\tau H_{k-1}(t) \, dt \, ,
\end{align*}
which implies for $|\beta|=2M$
\begin{align*}
   \sum_{0< \alpha\leq\beta} \frac{a_{\beta-\alpha}}{\alpha!}
\Big( \frac{x}{\sqrt{\cD}}\Big)^{\alpha}
&= \frac{1}{2^{2M} }\prod_{\beta_j>0} \frac{1}{\beta_j!}\Big( H_{\beta_j}
\Big( \frac{x_j}{\sqrt{\cD}}\Big)
-H_{\beta_j} (0)\Big)\\
&= \frac{1}{2^{2M} }\prod_{\beta_j>0}
\frac{2}{(\beta_j-1)!}\itg_0^{x_j/\sqrt{\cD}}H_{\beta_j-1}(t) \, dt \, .
\end{align*}
Consequently,
the remainder $\cR_{h,2M}$ takes the form
\begin{align*}
\Big(\frac{h \sqrt{\cD}}{ 2}\Big)^{N}(\pi \cD)^{-n/2} 
\sum_{|\beta|=2M}
   \sum_{m\in \R^n} {U_{\beta}(x,h m)}&\\
   \times{\ee}^{-|x-hm|^2/(h^2 \cD)}&
\prod_{\beta_j>0}
\frac{2}{(\beta_j-1)!}\itg_0^{z_j}H_{\beta_j-1}(t) \, dt \, ,
 \end{align*}
where we use the notation $z_j=(x_j-h m_j)/(h\sqrt{\cD})$.
Then Cramer's inequality for Hermite polynomials
\begin{align} \label{cramer}
|H_k(x)| \le  2^{k/2} \sqrt{k!} \, \ee^{\, x^2/2} 
\end{align}
(see \cite{SZ, IJ}),
leads to the estimate
\begin{align}
|\cR_{h,2M}(x)|& \le
\Big(\frac{h \sqrt{\cD}}{ 2}\Big)^{N}(\pi \cD)^{-n/2} \nonumber\\[-5pt]
& \quad \times
\sum_{|\beta|=2M}
   \sum_{m\in \R^n} |U_{\beta}(x,h m)|\, {\ee}^{-|x-hm|^2/(h^2 \cD)}
\prod_{\beta_j>0} \frac{2^{(\beta_j+1)/2}}{ \sqrt{(\beta_j-1)!} }
\itg_0^{|z_j|} \ee^{\, t^2/2} \, dt \nonumber \\
&\le \Big(\frac{h \sqrt{\cD}}{ 2}\Big)^{N}(\pi \cD)^{-n/2}
\sum_{|\beta|=2M} C_\beta
\sum_{m\in \R^n}|U_{\beta}(x,h m)| \,
S_\beta(x-h m) \, .\label{rm}
\end{align}
For the last inequality we use the notations
\begin{align*}
&S_\beta(x-h m)=\prod_{\beta_j=0}{\ee}^{-(x_j-hm_j)^2/(\cD h^2)}
\prod_{\beta_j>0} \frac{|x_j -hm_j|}{h \sqrt {2\cD}}{\ee}^{\, -(x_j
-hm_j)^2/(2h^2\cD)} \, , \\
&C_\beta = \prod_{\beta_j>0}
\frac{ 2^{\beta_j+1/2}}{  \sqrt{(\beta_j-1)!} } =
2^{2M}\prod_{\beta_j>0}
\sqrt{\frac{2}{(\beta_j-1)!} }
\end{align*}
and the estimate
\[
{\ee}^{\, -z_j^2}\itg_0^{|z_j|} \ee^{\, t^2/2} \, dt
\le |z_j|  \ee^{\,-z_j^2/2} \, .
\]
By \eqref{ualfa} we obtain for harmonic $u$
\begin{align*}
|U_{\beta}(x,h m)| &= \Big|\itg_{0}^{1}s^{N-1}\partial^{\beta} u(x+(1-s)(x-hm))\,
      ds \Big| \\
&\le \sum_{|\alpha|=0}^{\infty} \frac{\big|\de^{\beta+\alpha}u(x)\big|}{\alpha!}
\big|(x-hm)^{\alpha}\big| N \itg_{0}^{1}s^{N-1}(1-s)^{|\alpha|} ds \, ,
\end{align*}
which shows that
\begin{align*}
\sum_{m\in \R^n}|U_{\beta}(x,h m)| \,
S_\beta(x-h m)
&\le  \\ \sum_{|\alpha|=0}^{\infty}&
\frac{|\alpha|! \, N! }{ (N+|\alpha|)!} 
\frac{\big|\de^{\beta+\alpha}u(x)\big|}{ \alpha! }
\sum_{m\in \R^n}  S_\beta(x-h m) \big|(x-hm)^{\alpha}\big| .
\end{align*}
To get an upper bound of the last sum we write
\begin{align*}
&\sum_{m\in \R^n}  S_\beta(x-m)\big|(x-m)^{\alpha}\big|\\
&=\prod_{\beta_j=0} \sum_{m_j\in \R} |x_j -m_j|^{\alpha_j}{\ee}^{-(x_j -m_j)^2/\cD}
\prod_{\beta_j>0} \sum_{m_j\in \R}\frac{|x_j -m_j|^{\alpha_j+1}}{ \sqrt {2\cD}}
{\ee}^{\,  -(x_j  -m_j)^2/(2\cD)}\, ,
\end{align*}
and note that for $\cD > \cD_0$
\begin{align*}
& \sum_{m_j\in \R} |x_j -m_j|^{\alpha_j}{\ee}^{-(x_j -m_j)^2/\cD} \!\le c 
\cD^{(\alpha_j+1)/2}
\itg_0^\infty y^{\alpha_j} \ee^{\,-y^2} dy =
\frac{c\cD^{(\alpha_j+1)/2}}{2}\Gamma\Big(\frac{\alpha_j+1}{2}\Big),\\
& \sum_{m_j\in \R} \frac{|x_j -m_j|^{\alpha_j+1}}{ \sqrt {2\cD}}
{\ee}^{-(x_j -m_j)^2/(2\cD)} \le
\frac{c  \, (2\cD)^{(\alpha_j+1)/2}}{2}  \Gamma\Big(\frac{\alpha_j+2}{2}\Big).
\end{align*}
Hence we obtain the estimate
\begin{align*}
\sum_{m\in \R^n}  S_\beta(x-m)\big|(x-m)^{\alpha}\big|
\le \frac{c \,
\cD^{(|\alpha|+n)/2}}{2^n}
\prod_{\beta_j=0}\Gamma\Big(\frac{\alpha_j+1}{2}\Big) \prod_{\beta_j>0}
2 ^{(\alpha_j+1)/2} \, \Gamma\Big(\frac{\alpha_j+2}{2}\Big)
\end{align*}
with a constant $c$ independent of $\alpha$, $\beta$, and $\cD$. Now we use that
for $\alpha_j \ge 2$
\begin{align*}
\Gamma\Big(\frac{\alpha_j+1}{2}\Big)
\le \frac{\sqrt[4]{\pi}}{2^{\alpha_j/2}}  \sqrt{\alpha_j !}  \, ,
\end{align*}
which leads to
\begin{align*}
\sum_{m\in \R^n}  S_\beta(x-m)\big|(x-m)^{\alpha}\big|
\le
  c_1 \, \cD^{(|\alpha|+n)/2}
\prod_{\beta_j=0}\frac{ \sqrt{\alpha_j !}}{2^{\alpha_j/2}}
\prod_{\beta_j>0}
\sqrt{ (\alpha_j+1) !}
\end{align*}
with another constant $c_1$.
Thus we derive  from \eqref{rm}
\begin{align*}
&|\cR_{h,2M}(x)| \\
&\le C  (h \sqrt{\cD})^{N}\!\!\!
\sum_{|\beta|=2M}  \sum_{|\alpha|=0}^{\infty}
\frac{|\alpha|! \, |\beta|! \,\cD^{(|\alpha|+n)/2}}{ |\beta+\alpha|!} 
\frac{\big|\de^{\beta+\alpha}u(x)\big|}{ \alpha!}
\prod_{\beta_j=0}\frac{ \sqrt{\alpha_j !}}{2^{\alpha_j/2}}
\prod_{\beta_j>0}
\frac{ \sqrt{2 (\alpha_j+1) !}}{\sqrt{(\beta_j-1)!}}\\
&= C  (h \sqrt{\cD})^{N}\!\!\!
\sum_{|\beta|=2M}  \sum_{|\alpha|=0}^{\infty}
\frac{|\alpha|! \, |\beta|! \,\cD^{(|\alpha|+n)/2}}{ |\beta+\alpha|!}
\frac{\big|\de^{\beta+\alpha}u(x)\big|}{  \sqrt{\alpha ! \,   \beta!}}
\!\!\!\prod_{\beta_j=0} 2^{-\alpha_j/2}
\prod_{\beta_j>0} \sqrt{2 (\alpha_j+1)\beta_j}
\end{align*}
with a constant $C$ not depending on $M$, $\cD$, $h$, and $u$.
Since
\[
\frac{|\alpha|! \, |\beta|! }{ |\beta+\alpha|!} \sqrt{\frac{(\alpha+ \beta)!}{\alpha! \, \beta!}}
\le \frac{|\alpha|! \, |\beta|! }{ |\beta+\alpha|!} {\frac{(\alpha+ \beta)!}{\alpha! \, \beta!}}
\le 1 \, ,
\]
the remainder can be estimated by
\begin{align*}
&|\cR_{h,2M}(x)|
\le C  (h \sqrt{\cD})^{N}
\sum_{|\beta|=2M}  \sum_{|\alpha|=0}^{\infty}
  \frac{\big|\de^{\beta+\alpha}u(x)\big|}{\sqrt{(\alpha+\beta)!}}\cD^{(|\alpha|+n)/2}
\prod_{\beta_j>0} \sqrt{2 (\alpha_j+1)\beta_j}  \, .
\end{align*}
Thus, $|\cR_{h,2M}(x)| \to 0$ for any fixed $\cD$ if  \eqref{asssqu}
holds. 

\begin{remark}
{\em
The assertion of Theorem { \ref{thmarm}}
is a concrete realization of a general approximation result
for analytic functions.
Let $u$ be an entire function in $\C^n$ of order less than 2. 
Theorem 7.1 in \cite{MS2} 
states
that  the semi-discrete convolution
    \begin{equation*}
u_h(x) = \sum_{m \in {\Z}^n} u_{m} \, 
\ee^{- {|x-hm|^2}/(h^2\cD) }
     \end{equation*}           
with coefficients
     \begin{equation} \label{equ5e2}
   u_{m} := \itg_{\R^n} \ee^{ -  \pi^2 {\cD}|y|^2}
    \, u(h m + i \pi {\cD} h y) \, d y 
      \end{equation}           
differs from  $u$ by
\[
u_h(x)-u(x)=
\sum_{m \in \Z^{n}\setminus\{0\}}  \tilde u(x+i \pi h \cD m)
\ee^{ \,- \pi^2 \cD|m|^2} \ee^{ \,2 \pi i \langle   m,x\rangle/h}
  \, ,
\]
(cf. also \cite[Lemma~2.1]{MS5}).
It can be easily seen from \eqref{anext} and \eqref{equ5e2}
that the coefficients $u_{m}= (\pi \cD)^{-n/2} u(hm)$ if 
the restriction of $u$ to $\R^n$ is harmonic.
}
\end{remark}

\begin{figure}\begin{center}
         \includegraphics[width=2.4in,height=1.6in]{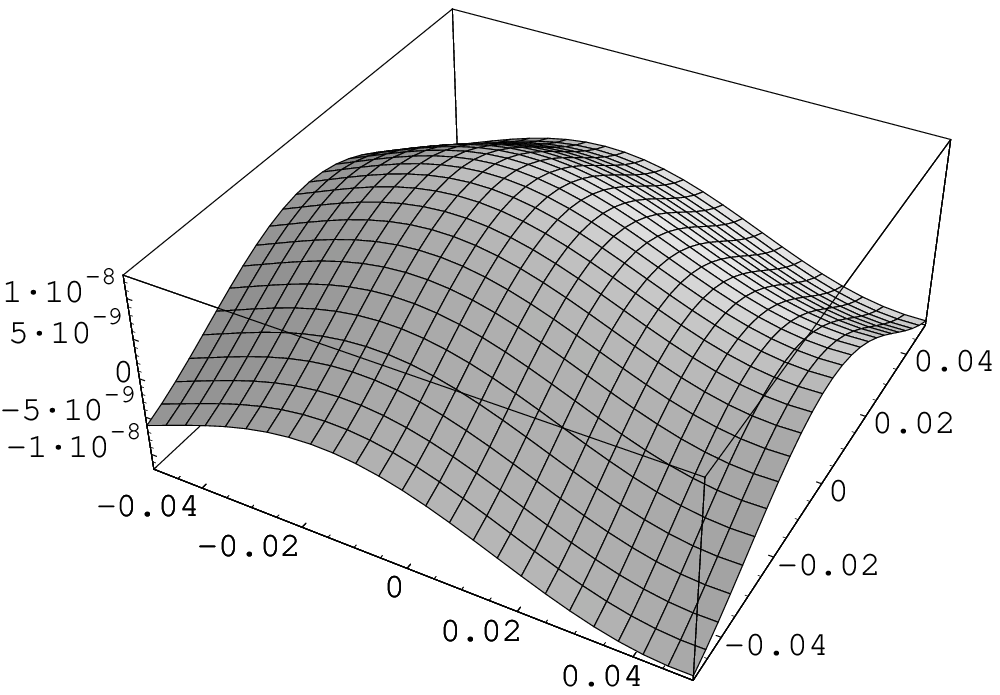}
      \includegraphics[width=2.4in,height=1.6in]{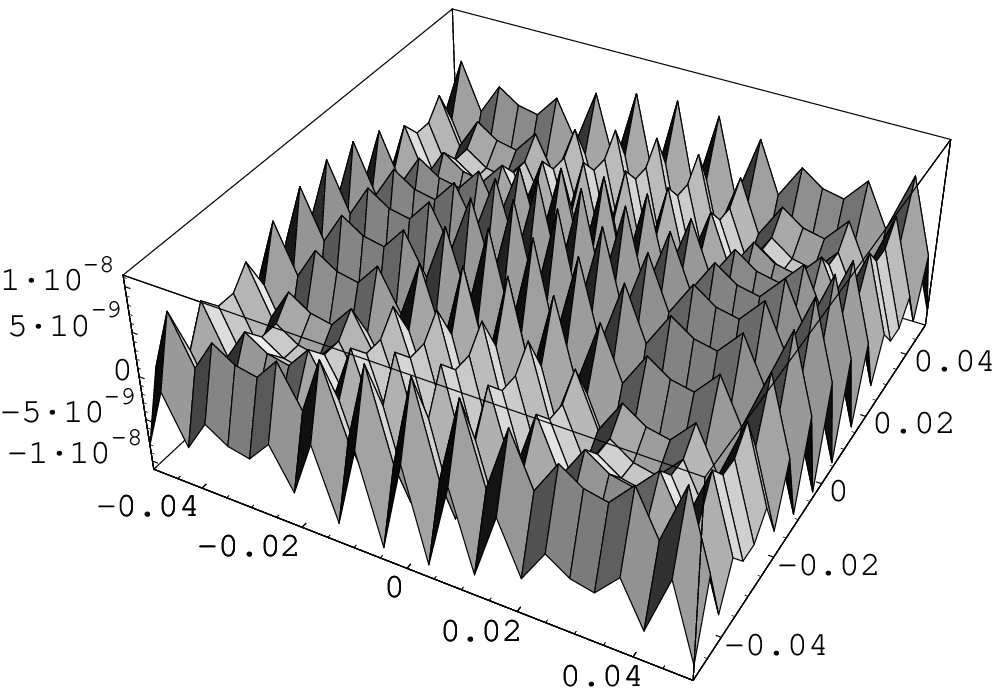}
\caption{ The graphs of $Mu-u$ with
$u(x_{1},x_{2})={\rm e}^{x_{1}}\cos x_{2}$,
$\cD=2$, $h=2^{-3}$ (on the left) and $h=2^{-7}$ (on the right).}
\label{harm1}
\end{center}
\end{figure}

\begin{figure}\begin{center}
     \includegraphics[width=2.4in,height=1.6in]{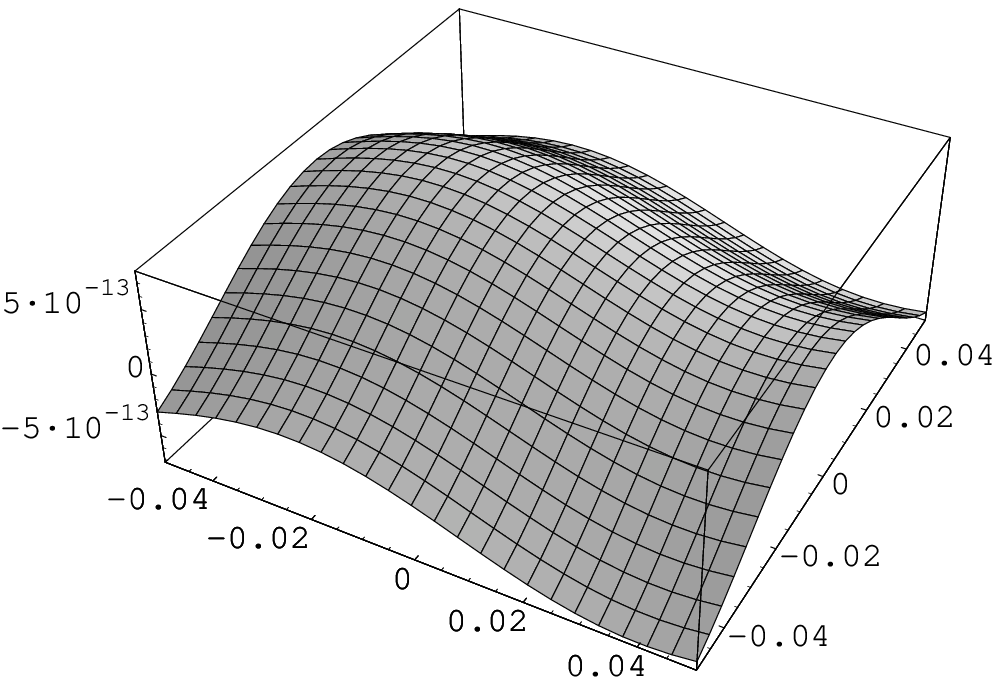}
      \includegraphics[width=2.4in,height=1.6in]{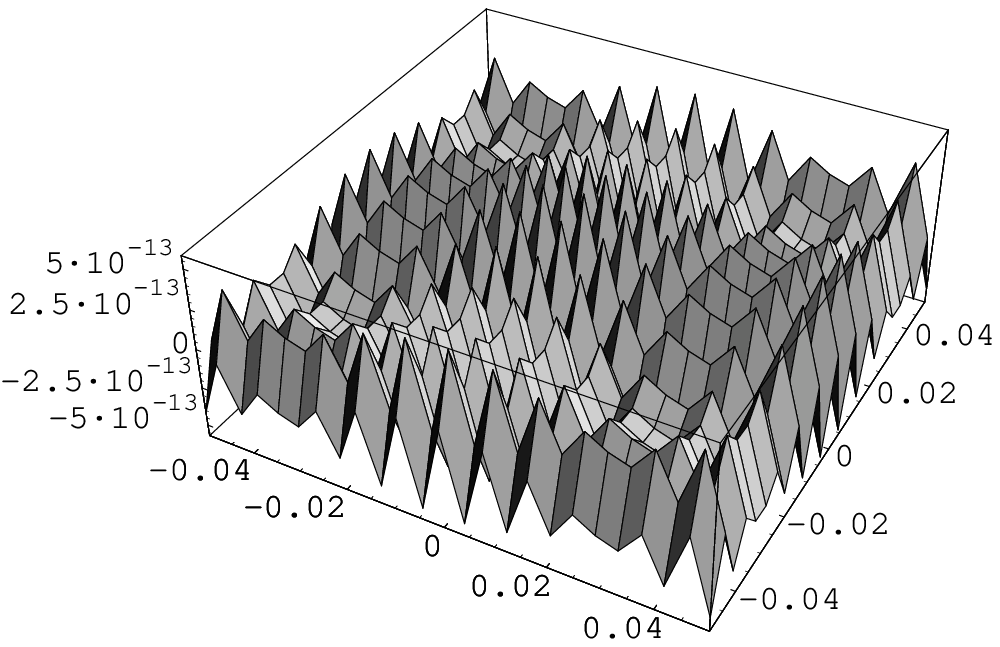}
\caption{The graphs of $Mu-u$ with $u(x_{1},x_{2})={\rm e}^{x_{1}}\cos
x_{2}$,
$\cD=3$, $h=2^{-3}$ (on the left) and $h=2^{-7}$ (on the right).}
\label{harm2}
\end{center}
\end{figure}

\begin{figure}\begin{center}
     \includegraphics[width=2.4in,height=1.6in]{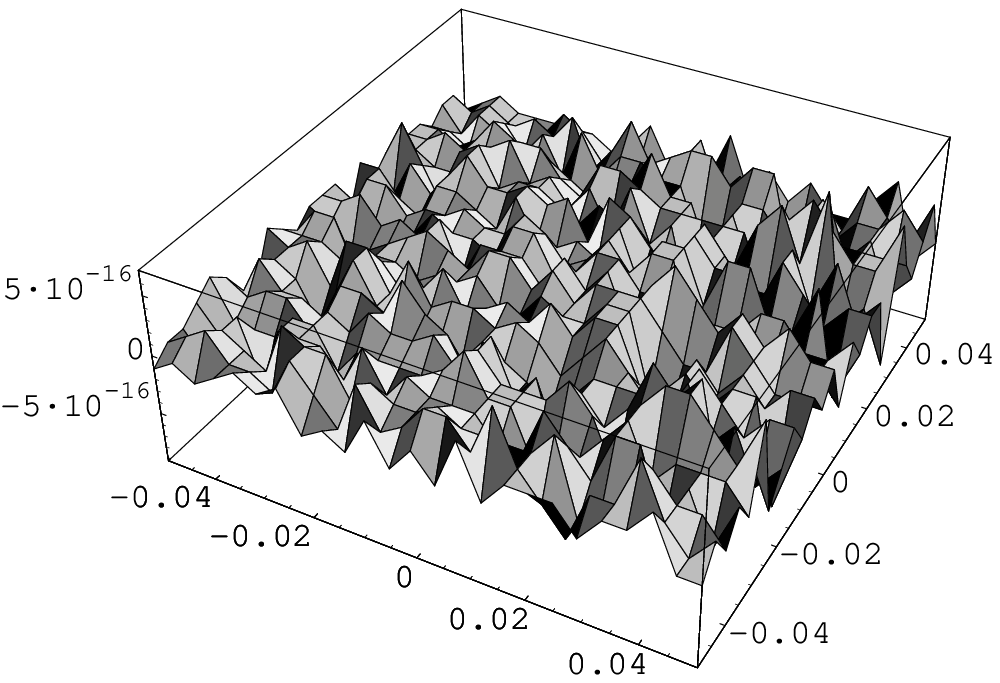}
      \includegraphics[width=2.4in,height=1.6in]{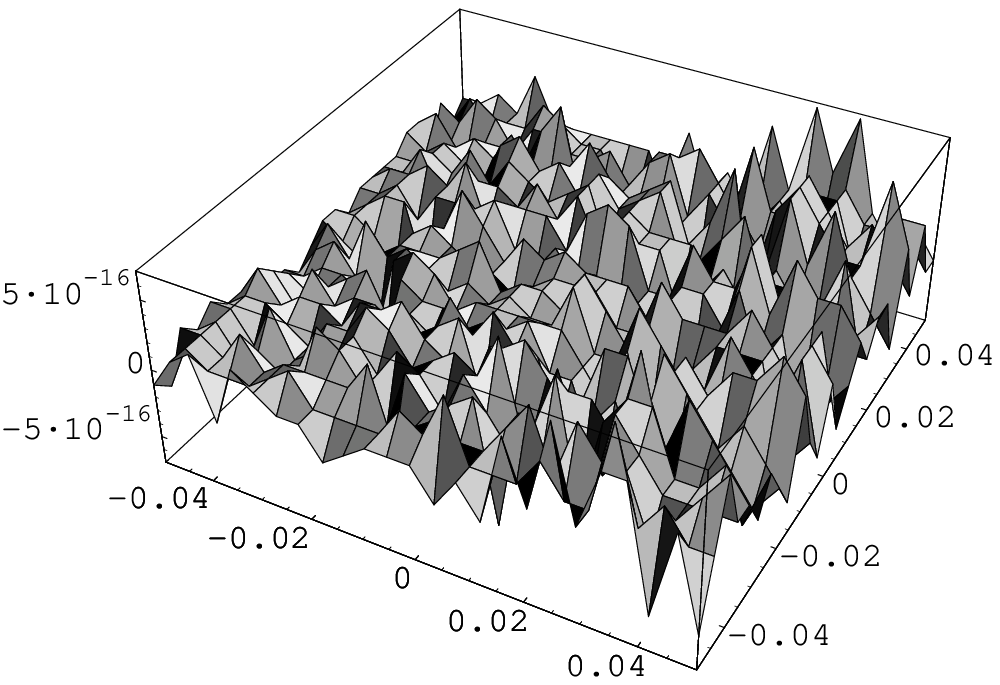}
\caption{The graphs of $Mu-u$ with $u(x_{1},x_{2})={\rm e}^{x_{1}}\cos
x_{2}$,
$\cD=4$, $h=2^{-3}$ (on the left) and $h=2^{-7}$ (on the right).}
\label{harm3}
\end{center}
\end{figure}

We have applied the simple quasi-interpolant \eqref{harmonic} 
to the harmonic function $u(x_{1},x_{2})={\rm e}^{x_{1}}\cos x_{2}$ in $\R^{2}$
by assuming $\cD=2$ (see Figure~\ref{harm1}), $\cD=3$  
(Figure~\ref{harm2}), $\cD=4$ (Figure~\ref{harm3}),
$h=2^{-3}$ and $2^{-7}$. The experiments confirm that the
quasi-interpolation error $Mu-u$ has reached its saturation
bound also  for large $h$ because it does not decrease if $h$ becomes
smaller.

Let now $u$ be harmonic  in some convex domain $\Omega \subset \R^n$ and we
consider the approximant
\begin{align} \label{harmonicdom}
M u(x)=(\pi \cD)^{-n/2}
\sum_{hm\in \Omega} u(hm)\, {\ee}^{-|x-h m|^2/(h^2 \cD)} \, .
\end{align}
\begin{theorem} \label{harmdom}
Suppose that the function $u$ is harmonic in a convex domain $\Omega
\subset \R^n$
and satisfies for a given $N=2M$
\[
c_u = \sum_{|\beta|=2M} \|\partial^\beta u \|_{L_\infty(\Omega)}
\prod_{\beta_j>0}
\sqrt{\frac{2}{(\beta_j-1)!}}
< \infty \, .
\]
Then for any $\e>0$ and subdomain $\Omega' \subsetneq \Omega$
there exists $\cD>0$ and $h>0$ such that
the quasi-interpolant \eqref{harmonicdom} provides
for all  $x \in \Omega'$ the estimate
      \begin{equation}\label{5.15}
          |u(x)- M u(x) |\leq C (h \sqrt{\cD})^{N} c_u
+\e \sum_{|\beta|=0}^{N-1} (h
          \sqrt{\cD})^{|\beta|}|\de^{\beta} u(x)| \, ,
      \end{equation}
where the constant $C$ depends only on the space dimension.
\end{theorem}
{\bf Proof.}
Analogously to the case $\Omega = \R^n$ we obtain
\[
|u(x)- M u(x) |\leq |\cR_{h,2M}(x)| +|\cE_{h,2M}(x)|
\]
with
\begin{align*}
|\cE_{h,2M}(x)|  &\le \,
|u(x)| \Big|(\pi \cD)^{-n/2}  \sum_{hm\in \Omega}
{\ee}^{-|x-hm|^2/(h^2\cD)} -1\Big| \\
+(\pi \cD)^{-n/2}& \sum_{|\beta|=1}^{2M-1}
\Big(\frac{h \sqrt{\cD}}{ 2}\Big)^{|\beta|}
   \frac{|\de^{\beta}u(x)|}{\beta! }
   \Big|\sum_{hm \in  \Omega} H_\beta
\Big( \frac{x-hm}{h \sqrt{\cD}}\Big)\, {\ee}^{-|x-hm|^2/(h^2\cD)}\Big| .
\end{align*}
and
\begin{align*}
|\cR_{h,2M}&(x)|
\le \Big(\frac{h \sqrt{\cD}}{ 2}\Big)^{N}(\pi \cD)^{-n/2}
\sum_{|\beta|=2M} C_\beta
\sum_{hm\in \Omega}|U_{\beta}(x,h m)| \,
S_\beta(x-h m) \, ,
\end{align*}
see \eqref{rm}. Since
\[
|U_{\beta}(x,h m)| \le \|\partial^\beta u \|_{L_\infty(\Omega)} \quad
\mbox{and} \quad  \sum_{hm\in \Omega} S_\beta(x-h m) \le c \,\cD^{n/2}
\]
with a constant $c$ depending only on $n$,
we get the inequality
\[
|\cR_{h,2M}(x)| \le C \, (h \sqrt{\cD})^{N} c_u \, .
\]
To estimate $|\cE_{h,2M}(x)|$ we use the functions $\sigma_\beta$ given by 
\eqref{defsigma}
and write
\begin{equation*}
\begin{split}
&(\pi \cD)^{-n/2}  \sum_{hm\in \Omega}
{\ee}^{-|x-hm|^2/(h^2\cD)} -1
=\sigma_0 \Big(\frac{x}{h}, \cD\Big)-(\pi \cD)^{-n/2}  \sum_{hm\notin \Omega}
{\ee}^{-|x-hm|^2/(h^2\cD)} \, , \\
&(\pi \cD)^{-n/2} \sum_{hm \in  \Omega} H_\beta
\Big( \frac{x-hm}{h \sqrt{\cD}}\Big)\, {\ee}^{-|x-hm|^2/(h^2\cD)} \\
&\hskip3cm = \sigma_\beta \Big(\frac{x}{h}, \cD\Big) -
   (\pi \cD)^{-n/2} \sum_{hm \notin  \Omega} H_\beta
\Big( \frac{x-hm}{h \sqrt{\cD}}\Big)\, {\ee}^{-|x-hm|^2/(h^2\cD)}\, .
\end{split}
\end{equation*}
Furthermore,
for $x \in \Omega$
we derive
\[
\bigg| (\pi \cD)^{-n/2} \sum_{hm \notin  \Omega} H_\beta
\Big( \frac{x-hm}{h \sqrt{\cD}}\Big)\, {\ee}^{-|x-hm|^2/(h^2\cD)}\bigg|
\le \delta_\beta(h^{-1} \dist(x,\partial \Omega), \cD) \, ,
\]
where  $\delta_\beta(r, \cD)$ , $r \ge 0$,
denotes the rapidly decaying
function
\[
\delta_\beta(r, \cD) =  \sup_{x\in\R^{n}}  \, (\pi \cD)^{-n/2}
\sum_{\substack{m \in \Z^n \\|x-m| > r}} \Big|H_\beta
\Big( \frac{x-m}{ \sqrt{\cD}}\Big)\Big|\, {\ee}^{-|x-m|^2/\cD} \, .
\]
Thus, for any
domain $\Omega$, fixed
parameter $h$ and multiindex $\beta$ we can find  a subdomain
$\Omega'_{\beta,h} \subsetneq \Omega$
such that
\begin{align} \label{xxx}
\sup_{x \in \Omega'_{\beta,h}} \delta_\beta(h^{-1} \dist(x,\partial \Omega), \cD) \le
\|\sigma_\beta (\cdot, \cD)\|_{L_\infty}\, ,
\end{align}
which gives
\begin{align*}
|\cE_{h,2M}(x)| \le & \,
2 \sum_{|\beta|=0}^{2M-1}
\Big(\frac{h \sqrt{\cD}}{ 2}\Big)^{|\beta|}
   \frac{|\de^{\beta}u(x)|}{\beta! }\|\sigma_\beta (\cdot, \cD)\|_{L_\infty}
\quad \mbox{for all} \quad x \in \bigcap_{|\beta|=0}^{2M-1} \Omega'_{\beta,h} \, .
\end{align*}
Now we have to choose $h$ such that
$\Omega' \subset \bigcap \Omega'_{\beta,h}$, which
is possible since $\Omega'_{\beta,h} \to \Omega$ as $h \to 0$.

\begin{remark}
{\em
Obviously the assertion of Theorems {\em \ref{thmarm}} and {\em \ref{harmdom}} can be extended
to the case that a solution $u$ of the second order equation
\[
\sum_{i,k=1}^{n} b_{ik} \de_{i}\de_{k} u =0
\]
in $\Omega$ is approximated by the quasi-interpolant
\begin{align*}
       M u(x)=\frac{(\det B)^{-1/2}}{(\pi \cD)^{n/2}}
\sum_{hm \in \Omega}u(hm)\,
\ee^{\, -\langle B^{-1}(x-hm),x-hm\rangle/(h^2\cD)}
\end{align*}
with the matrix $B=\{b_{ik}\}$.}
\end{remark}

\section{Approximation of derivatives}\label{deriv}

Here we study the approximation of derivatives
using Hermite quasi-interpolation operator \eqref{3.0}.
We introduce the continuous convolution (see \cite{MSbook})
\begin{equation}\label{conv}
C_{\delta}v(x)=\delta^{-n}\int_{\R^{n}}\cH\left(\frac{x-y}{\delta}\right) \cQ(-\delta
\,\de)v(y)dy
\end{equation}
where $\cQ(t)$ is the polynomial in \eqref{polinomio}.
\begin{theorem}\label{}
      Suppose that $\cH$ satisfies the decay condition \eqref{decaybis}
      with $K>L+n$, $L\in \N$, $L\geq N$.
For any $\e>0$ there exists $\cD>0$ such that for any
function $u\in W^{L}_{\infty}(\R^{n})$
\begin{equation}
\label{MC}
      |M u(x)-C_{h\dd} \,u(x)|\leq
      \e\,\sum_{|\gamma|=0}^{L-1}(h\dd)^{|\gamma|}|\de^{\gamma}u(x)|
    +c_{1}\,(h\dd)^{L}\sum_{{|\gamma|=L}}||\de^{\gamma}u||_{L_{\infty}} \, ,
\end{equation}
where the constant $c_{1}$ does not depend on $u$, $h$ and $\dd$.
\end{theorem}

{\bf Proof.}
Suppose that the function $u\in W^{L}_{\infty}(\R^{n})$.
The Taylor expansion
\eqref{3.2} with $N$ replaced by $L$
gives the following form of the quasi-interpolant $Mu$ in \eqref{3.1}
\begin{equation*}\begin{split}
       &M u(x)= \sum_{|\gamma|=0}^{N-1}a_{\gamma}\sum_{|\alpha|=0}^{L-1-|\gamma|}
       \frac{(- h \sqrt{\cD})^{|\alpha+\gamma|}}{\alpha!}\de^{\alpha+\gamma}u(x)
       \sigma_{\alpha}(\frac{x}{h},\cD,\cH) \\
&+ (-h \sqrt{\cD})^{L}
\sum_{|\gamma|=0}^{N-1}a_{\gamma}\!\!\!\sum_{|\alpha|=L-|\gamma|}  \frac{\cD^{-n/2}}{\alpha!}
    \!\!\sum_{m\in \Z^{n}} U_{\alpha+\gamma}(x,h m)
\left( \frac{x-h m}{h \sqrt{\cD}}\right)^{\alpha}\!\!\!
    \cH\left( \frac{x-h m}{h \sqrt{\cD}}\right).
\end{split}
\end{equation*}
Similarly the Taylor expansion of $u$ around $y$  leads to
\begin{equation}\label{conv2}
     \begin{split}
       C_{\delta}u(x)&= \sum_{|\gamma|=0}^{N-1}a_{\gamma}\sum_{|\alpha|=0}^{L-1-|\gamma|}
       (-\delta)^{|\alpha+\gamma|}\de^{\alpha+\gamma}u(x)
        \frac{ 1}{\alpha!}
       \int_{\R^{n}}\tau^{\alpha}\cH(\tau)d\tau\\
&+ (-\delta)^{L}
\sum_{|\gamma|=0}^{N-1}a_{\gamma}\sum_{|\alpha|=L-|\gamma|}  \frac{ 1}{\alpha!}
\int_{\R^{n}}\tau^{\alpha}\cH(\tau)U_{\alpha+\gamma}(x,x-\tau\delta)d\tau.
\end{split}
\end{equation}
     Setting $\delta=h \dd$,
we obtain
     \begin{equation*}
      \begin{split}
& Mu(x)-C_{h \dd}u(x)=
   \sum_{|\gamma|=0}^{N-1}a_{\gamma}\sum_{|\alpha|=0}^{L-1-|\gamma|}
   \frac{(- h \sqrt{\cD})^{|\alpha+\gamma|}}{\alpha!} \de^{\alpha+\gamma}u(x)
   \cE_{\alpha}(\frac{x}{h},\cD,\cH)\\
   & +(-h \sqrt{\cD})^{L}
\sum_{|\gamma|=0}^{N-1}a_{\gamma}\sum_{|\alpha|=L-|\gamma|}  \frac{ 1}{\alpha!}
   [ \cD^{-n/2}\!\!\!\sum_{m\in \Z^{n}} U_{\alpha+\gamma}(x,h m)
\left( \frac{x-h m}{h \sqrt{\cD}}\right)^{\alpha}\!\!\!
    \cH\left( \frac{x-h m}{h
   \sqrt{\cD}}\right)\\&-\itg_{\R^{n}}
    \tau^{\alpha}\cH(\tau)U_{\alpha+\gamma}(x,x-\tau\delta)d\tau].
       \end{split}
       \end{equation*}
       Then
\begin{equation*}\begin{split}
    & |Mu(x) \, -C_{h \dd}u(x)|\leq
   \sum_{|\gamma|=0}^{N-1}|a_{\gamma}|\sum_{|\alpha|=0}^{L-1-|\gamma|}
   \frac{( h \sqrt{\cD})^{|\alpha+\gamma|}}{\alpha!} |\de^{\alpha+\gamma}u(x)|
     \,\, |\cE_{\alpha}(\frac{x}{h},\cD,\cH)| \\
     &+(h \sqrt{\cD})^{L}
\sum_{|\gamma|=0}^{N-1}|a_{\gamma}|\sum_{|\alpha|=L-|\gamma|}\frac{
||\de^{\alpha+\gamma}u||_{L_{\infty}}}{\alpha!}
   \left(||\rho_{\alpha}(\cdot,\cD,\cH)||_{L_{\infty}}+
   \int_{\R^{n}}|\tau^{\alpha}\cH(\tau)|d\tau\right).
\end{split}
\end{equation*}

Proceeding as in the proof of Theorem \ref{UNO} we deduce
      \eqref{MC}. 

      \begin{theorem} \label{Cvv} If $\cH$ satisfies the conditions
\eqref{decaybis} with $K>N+n$ and
          \eqref{3.10},  then for any
          $v\in W^{N}_{\infty}(\R^{n})$
          \begin{equation}\label{conv3}
      |C_{\delta}v(x)-v(x)|\leq  c_{2} \delta^{N} \sum_{|\alpha|=N}\Vert
     \de^{\alpha}v\Vert_{L_{\infty}}.
\end{equation}
          \end{theorem}

{\bf Proof.} The representation \eqref{conv2} with $L=N$ gives
\begin{equation}\label{conv6}
      C_{\delta}v(x)=
      \sum_{|\alpha|=0}^{N-1}(-\delta)^{|\alpha|}\de^{\alpha}v(x)\sum_{\gamma\leq
     \alpha}\frac{a_{\alpha-\gamma}}{\gamma!}\int_{\R^{n}}\tau^{\gamma}\cH(\tau)d\tau+
      R_{\delta,\cH}(x)
\end{equation}
where the remainder $R_{\delta,\cH}$ satisfies
\begin{equation}\label{remconv}
     |R_{\delta,\cH}(x)|\leq \delta^{N}\sum_{|\alpha|=N}\Vert
     \de^{\alpha}v\Vert_{L_{\infty}} \sum_{\gamma\leq \alpha}\frac{|a_{\alpha-\gamma}|}{\gamma!}
     \int_{\R^{n}}|\tau^{\gamma}\cH(\tau)|d\tau\leq  c_{2} \delta^{N} \sum_{|\alpha|=N}\Vert
     \de^{\alpha}v\Vert_{L_{\infty}}.
\end{equation}
If condition (\ref{3.10}) holds, in view of (\ref{conv6}),
   we obtain \eqref{conv3}.

Theorem \ref{Cvv} leads immediately to the next corollary.

   \begin{corollary}\label{Cdudu}
       Suppose that $\cH$ satisfies conditions \eqref{decaybis} with $K>N+n$
       and
          \eqref{3.10}. Then
          for $u$ such that $\de^{\beta}u\in W^{N}_{\infty}(\R^{n})$,
   \begin{equation}\label{conv7}
       |C_{\delta}\de^{\beta}u(x)-\de^{\beta}u(x)|\leq
     c_{2}\delta^{N}\sum_{|\alpha|=N}\Vert
     \de^{\alpha+\beta}u\Vert_{L_{\infty}}.
\end{equation}
          \end{corollary}

If the derivative $\de^{\beta}\cH$
exists and satisfies the decay condition \eqref{decaybis} then
the convolution satisfies
          \begin{eqnarray}
      \label{CuC}
        \de^{\beta} C_{\delta}\,u(x)=(h
        \sqrt{\cD})^{-|\beta|}C_{\delta,\de^{\beta}\cH}\,u(x)=
        C_{\delta}\de^{\beta}u(x)
\end{eqnarray}
where
\[
C_{\delta,\de^{\beta}\cH}u(x)=\delta^{-n}\int_{\R^{n}}\de^{\beta}\cH\left(\frac{x-y}{\delta}\right) 
\cQ(-\delta
\,\de)u(y)dy
\]
and the quasi-interpolant $M u$ in (\ref{3.1})
satisfies  the equation
\begin{equation}\label{derivate}\de^{\beta}M u=(h
\sqrt{\cD})^{-|\beta|}M_{\de^{\beta}\cH}u
\end{equation}
with
\[ M_{\de^{\beta}\cH}u=
\cD^{-n/2} \sum_{m\in \Z^{n}}\de^{\beta}\cH\left( \frac{x-h m}{h
\sqrt{\cD}}\right)
\cQ\left(-h
\sqrt{\cD}\, {\partial }\right)  u(hm).
\]
Hence keeping in mind \eqref{CuC} and \eqref{derivate} we write the difference
\begin{align*}
(h &\sqrt{\cD})^{|\beta|}[\de^{\beta}M u(x)-\de^{\beta}u(x)]\\
&=    (h \sqrt{\cD})^{|\beta|}[\de^{\beta}Mu(x)-\de^{\beta}C_{\delta}u(x)]+
      (h \sqrt{\cD})^{|\beta|}[\de^{\beta}C_{\delta}u(x)-\de^{\beta}u(x)]\\
&=
      [M_{\de^{\beta}\cH}u(x)-C_{\delta,\de^{\beta}\cH}u(x)]+
       (h \sqrt{\cD})^{|\beta|}[C_{\delta}\de^{\beta}u(x)-\de^{\beta}u(x)].
\end{align*}
 From (\ref{conv7}) and (\ref{MC}) we obtain the following result.
\begin{theorem}
      Suppose that $\cH$ satisfies the conditions \eqref{decaybis} with
      $K>N+n$ and
      \eqref{3.10}.  Moreover suppose that the derivative $\de^{\beta} \cH$ exists
      and satisfies the conditions \eqref{decaybis} with $K>L+n$. The
      for any $\e>0$ there exits $\cD>0$ such that, for any $u\in
      W^{L}_{\infty}(\R^{n})$ with $L\geq N+|\beta|$,
\begin{align*}
    |\de^{\beta}M u(x)-\de^{\beta}u(x)|\leq& \, 
   \e\,\sum_{|\gamma|=0}^{L-1}(h\dd)^{|\gamma|-|\beta|}|\de^{\gamma}u(x)|\\
 +c_{1}\,&(h\dd)^{L-|\beta|}\sum_{{|\gamma|=L}}||\de^{\gamma}u||_{L_{\infty}}+
             c_{2}(h \dd)^{N}\sum_{|\alpha|=N}\Vert
\de^{\alpha+\beta}u\Vert_{L_{\infty}}.
\end{align*}
\end{theorem}

We deduce that formula (\ref{3.0}) gives the simultaneous
approximation of the derivatives $\partial^\beta u$ with the saturation
term $\e h^{-|\beta|}$.

\end{document}